\documentclass{AIMS}
\usepackage{amsmath}
  \usepackage{paralist}
  \usepackage{graphics} %% add this and next lines if pictures should be in esp format
  \usepackage{epsfig} %For pictures: screened artwork should be set up with an 85 or 100 line screen
\usepackage{graphicx}  \usepackage{epstopdf}%This is to transfer .eps figure to .pdf figure; please compile your paper using PDFLeTex or PDFTeXify.

\graphicspath{ {images/} }
 \usepackage[colorlinks=true]{hyperref}
\hypersetup{urlcolor=blue, citecolor=red}

\newtheorem{theorem}{Theorem}[section]

\newtheorem{conjecture}{Conjecture}

\theoremstyle{definition}
\newtheorem{definition}[theorem]{Definition}

\newcommand{\scalor}{0.8}
\newcommand{\scalorDif}{0.95}

\title[Spectral Asymptotics of the Laplacian] %Use the shortened version of the full title
      {Spectral Asymptotics of the Laplacian on Surfaces of Constant Curvature}

\author[Timothy Murray and Robert S Strichartz]{}


 \email{tsmurra2@illinois.edu}
 \email{str@math.cornell.org}

\thanks{$^*$ Research of the second author is supported in part by the National Science Foundation, Grant DMS - 1162045}

\begin{document}
\maketitle

% Enter the first author's name and address:
\centerline{\scshape Timothy Murray}
\medskip
{\footnotesize
% please put the address of the first author
 \centerline{Department of Industrial and Systems Engineering}
   \centerline{University of Illinois at Urbana-Champaign}
   \centerline{Transportation Building, Urbana, IL 61801, USA}
} % Do not forget to end the {\footnotesize by the sign }

\medskip

\centerline{\scshape Robert S Strichartz$^*$}
\medskip
{\footnotesize
 % please put the address of the second  and third author
 \centerline{Department of Mathematics}
   \centerline{Cornell University}
   \centerline{Malott Hall, Ithaca, NY 14853, USA}
}

\bigskip

% The name of the associate editor will be entered by an editorial staff
% "Communicated by the associate editor name" is not needed for special issue.
 \centerline{(Communicated by the associate editor name)}

%The abstract of your paper
\begin{abstract}
The purpose of this paper is to explore the asymptotics of the eigenvalue spectrum of the Laplacian on 2 dimensional spaces of constant curvature, giving strong experimental evidence for a conjecture of the second author \cite{strichartz2016}. We computed and analyzed the eigenvalue spectra of several different regions in Euclidean, Hyperbolic, and Spherical space under Dirichlet, Neumann, and mixed boundary conditions and in particular we found that the average of the difference between the eigenvalue counting function and a 3-term prediction has the expected nice behavior. All computational code and data is available on our companion website \cite{murray2015results}.
\end{abstract}

\section{Introduction}
The purpose of this paper is to present numerical evidence for a conjecture on the spectral asymptotics of the Laplacian on a surface of constant curvature presented by the second author in \cite{strichartz2016}. The conjecture was initially limited to surfaces of either zero or constant  positive curvature, but we present strong evidence that a version of it should also be valid in the case of constant negative curvature. Results of Bleher \cite{bleher1994} show that it cannot be valid for surfaces of variable curvature.

We consider surfaces $S$ of finite area $A$ with boundary $\partial S$ of finite perimeter $P$ 
that is made up of a finite number of smooth curves meeting at angles ${\{\theta_{j}\}}$. Simple
examples are triangles and discs in either the Euclidean plane, the sphere, or hyperbolic 2-space.
We will also look at more complicated examples where $S$ is not simply connected and has 
non-convex boundary. We consider the standard Laplacian $\Delta$ with either Dirichlet ($D$),
Neumann ($N$), or mixed boundary conditions ($D$ on a portion of the boundary with perimeter
${P_{D}}$, and $N$ on the remaining portion of the boundary with perimeter ${P_{N}}$). We let
${\{\lambda_{j}\}}$ be the set of eigenvalues ${-\Delta u_{j} = \lambda_{j}u_{j}}$ repeated 
according to multiplicity, so the ${\lambda_{j}}$ are all nonnegative and ${\lambda_{j}\rightarrow\infty}$
as ${j\rightarrow\infty}$. The eigenvalue counting function is defined as 
\begin{equation}
{N(t) = \# \{\lambda_{j} \leq t\}}\tag{1.1} \newline 
\end{equation}
(Note that some references will use $\lambda_{j}\leq t^2$ instead). The well-known Weyl asymptotic formula $N(t) \sim \frac{A}{4\pi}t$ was refined by Ivvii \cite{ivrii1984precise} to
\begin{equation}
\tag{1.2} {N(t) = \frac{A}{4\pi}t + \frac{P_{N}-P_{D}}{4\pi}t^{1/2}} + O(t^{1/2}) \newline
\end{equation}
In \cite{strichartz2016} we proposed a still more refined asymptotic
\begin{equation}
\tag{1.3} {\widetilde{N}}(t) = \frac{A}{4\pi}t + \frac{P_{N}-P_{D}}{4\pi}t^{1/2} + C \newline
\end{equation} 
where the constant $C$ will be explained in Definition 1.1 below. However, it is impossible to see the constant from $N(t)$ alone, since it is expected that 
\begin{equation}
\tag{1.4} D(t) = N(t) - {\widetilde{N}(t)}
\end{equation}
will have at least growth $O(t^{1/4})$. Instead we consider the ordinary average error
\begin{equation}
\tag{1.5} A(t) = \frac{1}{t}\int_0^t D(s)\,ds \newline
\end{equation} We conjecture that this is bounded and decays on the order of $O(t^{-1/4})$ in the flat or negative curvature case. See Conjecture 1.2 below for a more detailed description. For this to be valid we need the correct value for the constant. 

We note that a different kind of average, the trace of the heat kernel
\begin{equation}
\tag{1.6}h(t) = - \sum_j e^{-t\lambda_{j}} \newline
\end{equation} (as $t\rightarrow 0$) has been extensively studied, beginning with the famous paper of Mark Kac \cite{kac1996hear} and continuing with \cite{vandenBerg1990}, \cite{buser1992geometry}, \cite{gilkey2003asymptotic}, \cite{mckean_jr.1967}, and \cite{stewartson_waechter_1971}. A related "logarithmic Gaussian averaged error estimate" is studied by Brownell in \cite{brownell1957extended} (see also \cite{baltes1976spectra} for a discussion of this). These are smoother type averages than the ones we consider, and in particular they involve the entire spectrum. Because they are smoother averages, they effectively erase some of the interesting detail that the rougher averages see. It is straightforward to obtain the smooth average results from the rough average results, and in particular our formula (1.3) for the refined asymptotics is consistent with the earlier results. To go in the reverse direction requires using a Tauberian theorem that only yields the original Weyl asymptotics. The average $A(t)$ is a special case of Riesz means, which have also been studier, starting with H\"{o}rmander \cite{hörmander1968}. See also \cite{Fulling1999} for an extensive survey of this approach.

The method we use to numerically approximate the spectrum of the Laplacian is extremely straightforward. We use the finite element solver built into MATLAB. For surfaces in the plane we just have to give a description of the boundary. For
surfaces in the hyperbolic plane or sphere we us a conformally flat coordinate system so the surface Laplacian becomes a scalar multiple of the Euclidean Laplacian. By using the mesh refinement option we obtain better and better approximations of
smaller initial segments of the spectrum. Given the computation time constraints, this allows us to get confident approximations for only a couple hundred eigenvalues. We then use an ad hoc extrapolation method on the sequence of approximations with
increasing refinements to get a slightly improved final approximation. We were pleasantly surprised to see that this small peek at an initial segment of the spectrum already yields strong evidence for the conjecture. In other words, it appears that the asymptotic regimine kicks in very early in the game. In the case of the Euclidean disc we have a better alternative method, since there the eigenvalues are given explicitly as squares of zeroes of Bessel functions of the first kind (D) or zeroes of derivatives of Bessel functions of the first kind (N). This allows us to go higher up in the spectrum with greater accuracy, and serves as a check on the size of the error obtained by the cruder method. Another check on error size is provided by doing the computations for the few triangles where the exact spectrum is known. 

We now present the details concerning the constant $C$ in (1.3)

\begin{definition}
Let $K_{2}(S)$ denote the curvature of $S$, which is assumed to be constant, and let $K_{1}$ denote the curvature function on the smooth pieces of $\partial S$ as viewed from $S$.
\end{definition}
Further define
\begin{equation}
\tag{1.7} \varphi (\theta) = \frac{1}{24}\left(\frac{\pi}{\theta}-\frac{\theta}{\pi}\right) \newline 
\end{equation}Then
\begin{equation}
\tag{1.8}C = C_1 + C_2 + C_3 \newline 
\end{equation}
where
\begin{equation}
\tag{1.9}C_3 = \frac{1}{12\pi}A K_2(S) \\
\end{equation}
\begin{equation}
\tag{1.10}C_2 = \frac{1}{12\pi}\int_{\partial S} K_1 \\
\end{equation}
\begin{equation}
\tag{1.11}C_1 = \sum_j \varphi(\theta_j) \\
\end{equation}
in the case of $D$ or $N$ boundary conditions throughout, or
\begin{equation}
\tag{1.12}C_1 = \sum \varphi(\theta_j') + \sum (\varphi(2\theta_j'')-\varphi(\theta_j''))
\end{equation}
for mixed boundary conditions, where the corner angles are sorted into $\{\theta_j'\}$ where the same type of boundary condition is imposed on both sides of the corner, and $\{\theta_j''\}$ where opposite type boundary conditions are imposed on the two side arcs.\par We note that in [S] we also allowed a finite number of cone point singularities on $S$ with cone angles  $\{\alpha_j \}$, and these contributed an additional term
\begin{equation}
\tag{1.13}\sum 2\varphi(\alpha_j/2) \newline
\end{equation}
 to $C_1$. However, we are unable to do our computations if there are cone point sungularities, so we can't test the conjecture in such cases. 

\begin{conjecture}\label{conj_curv_neg}
Assume the curvature of $S$ is zero or negative. Then there exists a uniformly almost periodic function $g$ such that 
\begin{equation}
\tag{1.14}A(t) = g(t^{1/2})t^{-1/4}+O(t^{-1/2})
\end{equation}
 as $t\rightarrow\infty$. In the case of zero curvature the almost periodic function $g$ has mean value zero.
\end{conjecture} 

\begin{conjecture}\label{conj_curv_pos}
Assume the curvature of $S$ is positive. Then there exists a uniformly almost periodic funtion $g$ of mean value zero such that
\begin{equation}
\tag{1.15}A(t) = g(t^{1/2}) + O(t^{-1/2})
\end{equation}as $t\rightarrow 0$.
\end{conjecture}

We note that in \cite{strichartz2016} we conjectured that (1.14) and (1.15) are the first terms in an asymptotic expansion, but we are unable to test this here. Indeed, we cannot test the rate of decay in (1.14) and (1.15), since we don't know what $g$ should be. So basically we will observe that $t^{1/4}A(t)$ in case the case of conjecture \ref{conj_curv_neg} and $A(t)$ in the case conjecture \ref{conj_curv_pos} appear to be bounded functions of $t^2$ with mean value zero that could reasonably be almost periodic. Since almost periodicity is a global property, there is no way to test it by examining a small portion of the graph. We will observe, however, that there is no discernable difference between examples where the almost periodicity is known to be true, and all the other examples.

This paper is organized as follows: in section \ref{sec_test_ex} we perform our experimental computations for examples where the spectrum is known exactly, two Euclidean triangles and the Euclidean disc with Dirichlet and Neumann boundary conditions. We introduce the six part graphical display of data that will be used throughout the paper (except for the spherical surfaces in section \ref{sec_sph}). The reader will be able to see at a glance both confirmation of the predicted behavior and deviations due to computational error. In section \ref{sec_flat} we examine many examples of flat surfaces, including surfaces with mixed boundary conditions, surfaces that are not convex, and surfaces that are not simply connected. In section \ref{sec_hyp} we study hyperbolic surfaces, both triangles and discs. We see here experimental evidence that the conjecture for flat surfaces carries over into this case. In section \ref{sec_sph} we study spherical surfaces. Since the conjecture is different in this case (with no decay in $A(t)$) we use a five part graphical display. We give a discussion of all our results in section \ref{sec_discus}. We also mention the interesting question of the behavior of differences of consecutive eigenvalues. We have gathered data for all the examples studied here, and present a small selection of it. At this time we are not able to propose any conjectures.

The website \cite{murray2015results} contains the complete data on all the examples discussed here, as well as many other related examples. Additionally, a zip-file of all of our code is available for download. Automated scripts to generate each set of eigenvalues for an arbitrary number of refinements are available. However, please note that as much of each experiment was done through in-console manipulations there is no one unified script or function to generate the predicted eigenvalues or graphs once the initial refinements are performed.

\section{Some Test Examples}\label{sec_test_ex}
In this section we discuss our results for a few examples of surfaces where the spectrum is known exactly.

\subsection{Euclidean right isosceles triangle with Dirichlet boundary conditions}

The set of eigenvalues is $\pi^{2}(j^{2}+k^{2})$ for all pairs $(j,k)$ of distinct positive integers. We will normalize all eigenvalues by dividing by $\pi^2$ so that we are dealing with integer values. In Table 1 we show the data for the first 10 eigenvalues (the full table is on the website \cite{murray2015results}). In the first column we show the initial MATLAB computation of $\frac{\lambda_j}{\pi^2}$. In the next 6 columns we show the same value after successive refinements of the mesh. So the initial value of $\frac{\lambda_{10}}{\pi^2}$ is 44.931704, which is quite far from the true value of 37, but by the 6$^{th}$ refinement the approximation has improved to 37.001949. The next column is our predicted value obtained by fitting the data $x_{n}$ for refinements $n = 4,5,6$ to $x_{n} = x + cr^{n}$ and taking $x$ for the prediction. In this case the prediction is 37.000001. If we look further up in the spectrum we can see eigenvalues with multiplicity 2. For example $\frac{\lambda_{133}}{\pi^{2}} = \frac{\lambda_{134}}{\pi^{2}} = 377$. At refinement 4 the two values are 380.03292 and 380.1314. Not very close to each other and far off from the true value. At refinement 5 the two values are 377.7568 and 377.7816, closer to the true value but still not too close to each other. The predicted values are 376.9999 and 376.9998. Even though the order gets switched, the error is still quite acceptable. 

\begin{table}[ht]\label{tab_euc_right_iso_tri}
\centering
	\begin{tabular}{| l | l | l | l | l | l | l | l | }
	\hline
	 & Initial & 1 & 2 & 3 & 4 & 5 & Predicted \\ \hline
1 & 5.13589 & 5.03479 & 5.0088 & 5.00221 & 5.000134 & 5.000138 & 5 \\ \hline
2 & 10.5735 & 10.1448 & 10.0364 & 10.0091 & 10.0006 & 10.0006 & 10 \\ \hline
3 & 13.9042 & 13.2281 & 13.0573 & 13.0143 & 13.0009 & 13.0009 & 13 \\ \hline
4 & 18.6783 & 17.4194 & 17.1051 & 17.0263 & 17.0016 & 17.0016 & 17 \\ \hline
5 & 22.3425 & 20.5806 & 20.1451 & 20.0363 & 20.0022 & 20.0022 & 20 \\ \hline
6 & 28.5140 & 25.8760 & 25.2190 & 25.0548 & 25.0034 & 25.0034 & 25 \\ \hline
7 & 29.7992 & 26.9473 & 26.2370 & 26.0593 & 26.0037 & 26.0037 & 26 \\ \hline
8 & 33.5825 & 30.1526 & 29.2891 & 29.0724 & 29.0045 & 29.0045 & 29 \\ \hline
9 & 40.6045 & 35.6485 & 34.4114 & 34.1029 & 34.0064 & 34.0064 & 34 \\ \hline
10 & 44.9317 & 38.9934 & 37.4981 & 37.1246 & 37.0078 & 37.0077 & 37 \\ \hline
... & ... & ... & ... & ...& ... & ... & ...\\ \hline
133 & 0 & 0 & 0 & 389.656 & 377.756 & 377.756 & 376.9999 \\ \hline
134 & 0 & 0 & 0 & 398.873 & 377.781 & 377.781 & 376.9998 \\ \hline
	\end{tabular}

	\caption{Euclidean right isosceles triangle}

\end{table}

In Figure 1 we show the graphs of
\begin{enumerate}
\item $N(t)$
\item $D(t) = N(t) - \widetilde{N}(t)$
\item $A(t) = \frac{1}{t}\int_{0}^{t}D(s)ds$
\item $t^{\frac{1}{4}}A(t)$
\item $t^{\frac{1}{4}}A(t^{2})$
\item $\frac{1}{t-a}\int_{a}^{t}s^{\frac{1}{2}}A(s^{2})ds$ for $a = $ the highest predicted eigenvalue divided by 16, removing the first $\frac{1}{4}$ of graph 5 from figuring into graph 6 and eliminating potential early extreme values so that it converges to 0 more quickly.
\end{enumerate}
We will use this set of six graphs for each Euclidean and hyperbolic region which we analyze. The $x$-scales of the first four graphs were picked to use all predicted eigenvalues with an acceptable level of error, usually between the first 120 and 150 eigenvalues (the number used is in the $y$-axis of the first graph). When the true values are known, we use those and frequently use more than 150, as in figure 2 where the first 1000 eigenvalues are used. The scale of the $x$-axis in the fifth and sixth graphs is approximately the square root of the scale of the $x$-axis of the first four graphs. 

For Figure 1 We used the exact values for the first approximately 150 eigenvalues. A quick look at these graphs yields some simple observations. The graph 1 shows that $N(t)$ grows approximately linearly, while 2 shows that $D(t)$ grows at a relatively slow rate. the graph 3 suggests that $A(t)$ is converging to 0 at a slow rate, while graph 4 confirms that $O(t^{-\frac{1}{4}})$ is a plausible decay rate. The function in graph 5 is known to be converging to an almost periodic function g(t) (see \cite{strichartz2016}), but this is not apparent from the graph. Presumably the almost periods are too large to show up in the range of data we have plotted. On the other hand, graph 6 gives strong evidence that the almost periodic function has mean value 0.

\begin{figure}[ht]\label{count_euc_right_iso_tri}
\centering
\includegraphics[width = \scalor\textwidth]{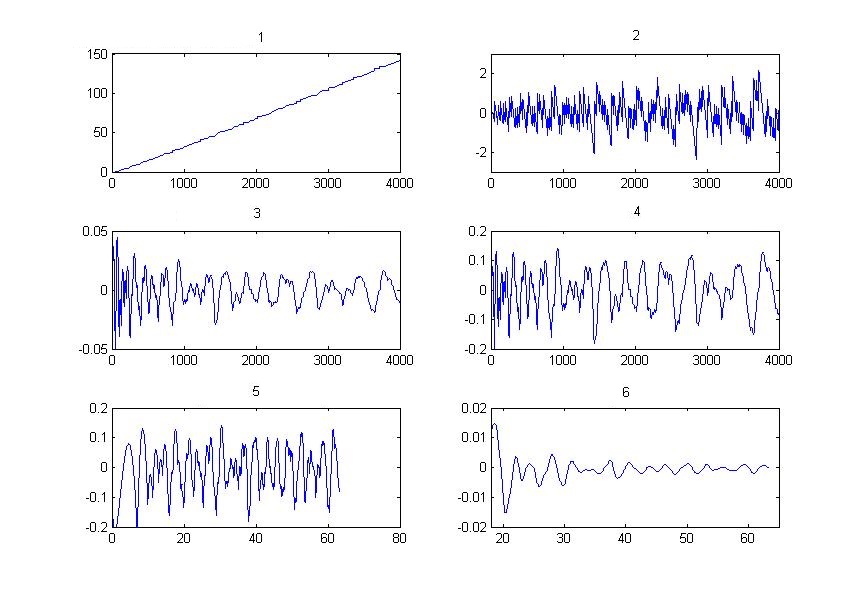}\\
\caption{Euclidean right isosceles triangle}       
\centering
\end{figure}

%In Figure 2.2 we display the set of differences $\frac{(\lambda_{k+1} - \lambda{k})}{\pi^{2}}$ , both as a cumulative sum and a %histogram. Since these differences are all integers there is no question about the choice of bin sizes for the histogram. It is not %clear from the limited data how the histogram will behave in the limit. The fact that the bin values are very far from being 
%decreasing is very intriguing. Ultimately it should be an interesting problem in number theory to find any pattern.
%\begin{center}
%\includegraphics[scale=.45]{EucDrightIsoTriDifferenceGraphs}
%Figure 2.2
%\end{center} 

\subsection{The Euclidean equilateral triangle with Dirichlet boundary conditions} 

Here the eigenvalues are known to be the values $(\frac{4}{3}\pi)^{2}(j^{2} + k^{2} +jk)$ for the positive integers $j,k$. This typically produces multiplicity 1 when $j = k$ and multiplicity 2 when $j \neq k$. Here we only used 5 refinements. Table 2 and Figure 2 show the same information for this example as before. For $\lambda_{119} = \lambda_{120} = 219$ our predicted values are 219.0009455 and 219.0005973 while on the $5^{th}$ refinement they are 219.4644686 and 219.4662058. The qualitative features of Figure 2 are much the same as that of Figure 1.

\begin{figure}[ht]\label{count_euc_equ_tri}
\centering
\includegraphics[width=\scalor\textwidth]{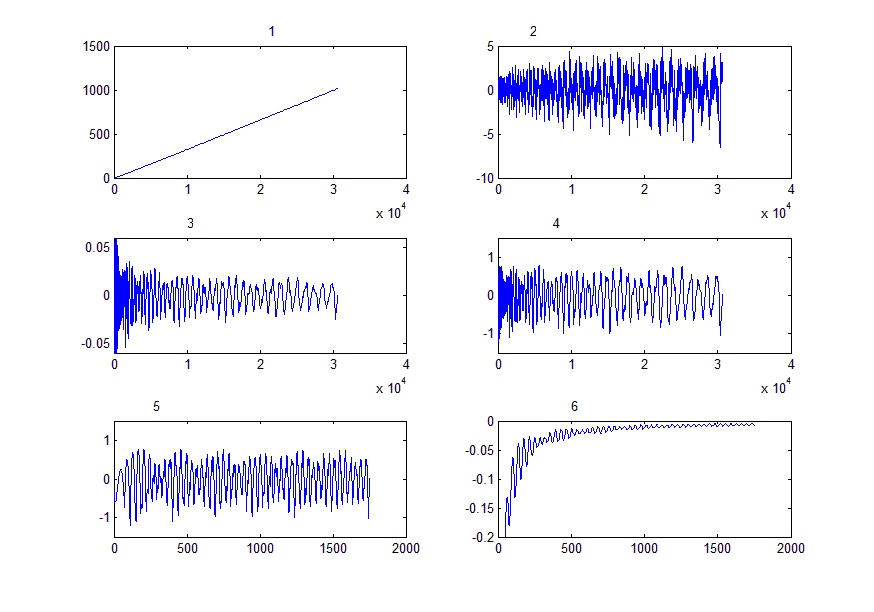}\\
\caption{Euclidean equilateral triangle}
\end{figure} 

%\begin{center}
%\includegraphics[scale=.6]{EucDEquTriDifferenceGraphs}
%Figure 2.4
%\end{center} 

\begin{table}[ht]\label{tab_euc_equ_tri}
	\centering
	\begin{tabular}{| l | l | l | l | l | l | l | l | }
	\hline
	 & Initial &  1 & 2 & 3 & 4 & 5 & Predicted \\ \hline
	1 & 3.08733 & 3.02224 & 3.005611 & 3.001408 & 3.000088 & 3.000088 & 3.000001 \\ \hline
2 & 7.44103 & 7.11170 & 7.028133 & 7.007054 & 7.000441 & 7.000441 & 7.000003 \\ \hline
3 & 7.49871 & 7.12432 & 7.03113 & 7.007791 & 7.000487 & 7.000487 & 7.000003 \\ \hline
4 & 13.4606 & 12.3613 & 12.09015 & 12.02254 & 12.00141 & 12.00141 & 12.00001 \\ \hline
5 & 14.6422 & 13.4084 & 13.10213 & 13.02555 & 13.0016 & 13.0016 & 13.00001 \\ \hline
6 & 14.7099 & 13.4240 & 13.10604 & 13.02654 & 13.00166 & 13.00166 & 13.00001 \\ \hline
7 & 22.3866 & 19.8479 & 19.2122 & 19.05312 & 19.00332 & 19.00332 & 19.00002 \\ \hline
8 & 22.7581 & 19.9425 & 19.23499 & 19.05874 & 19.00367 & 19.00367 & 19.00002 \\ \hline
9 & 25.3291 & 22.0593 & 21.26293 & 21.06566 & 21.0041 & 21.0041 & 21.00002 \\ \hline
10 & 25.3551 & 22.0705 & 21.26635 & 21.06656 & 21.00416 & 21.00416 & 21.00002 \\ \hline

	\end{tabular}
	\caption{Euclidean equilateral triangle}
\end{table}

\subsection{Euclidean disc with Dirichlet boundary conditions}

We take the radius to be one since all discs have eigenvalues that scale by the radius. In this case the eigenvalues are the squares of the zeroes of the Bessel functions $J_{k}$ for nonnegative integers $k$ with multiplicity one for $k = 0$ and multiplicity two for $k \geq 1$. It is possible to get accurate values of these zeros so we have exact values for the first 660 eigenvalues. In this example $\widetilde{N}(t) = \frac{1}{4}t - \frac{1}{2}\sqrt{t} + \frac{1}{6}$. 

% Note, however, that the validity of the conjecture is not known in this example. In Table 2.3 we show the first 10 computed approximations with 4 refinements. The table also lists the zeros of the Bessel function and the values $k$ and $n$ for $J_{k}(z_{n}) = 0$. The accuracy of our predictions are worse than in the triangle examples. For example $\lambda_{149} = \lambda_{150} = 646.7022512 = (z_{6})^{2}$ for $J_{5}(z_{6}) = 0$ has fairly accurate predicted values 646.7227695 and 646.7146116, while $\lambda_{212} = \lambda_{213} = 910.7757336 = (z_{5})^{2}$ for $J_{11}(z_{5}) = 0$ has unacceptable predicted values 910.4406915 and 906.2756139.

Figure 3 displays the same graphs as before using the exact values. We note that graph 5 is just as plausibly an asymptotic almost periodic function as the same graphs in the triangle cases where we know the function is asymptotically almost periodic. On the other hand, graph 6 shows a much slower rate of decay than in the triangle cases. It is still plausible that this gives supportive evidence that the presumed almost periodic function has mean value zero, but the evidence is not decisive. 

%Figure 2.6 show the eigenvalue difference. In this case the histogram seems to have a different appearance than for the triangle cases, but now the result is sensitive to the choice of bin size.

\begin{figure}[ht]\label{count_euc_disc}
\centering
\includegraphics[width=\scalor\textwidth]{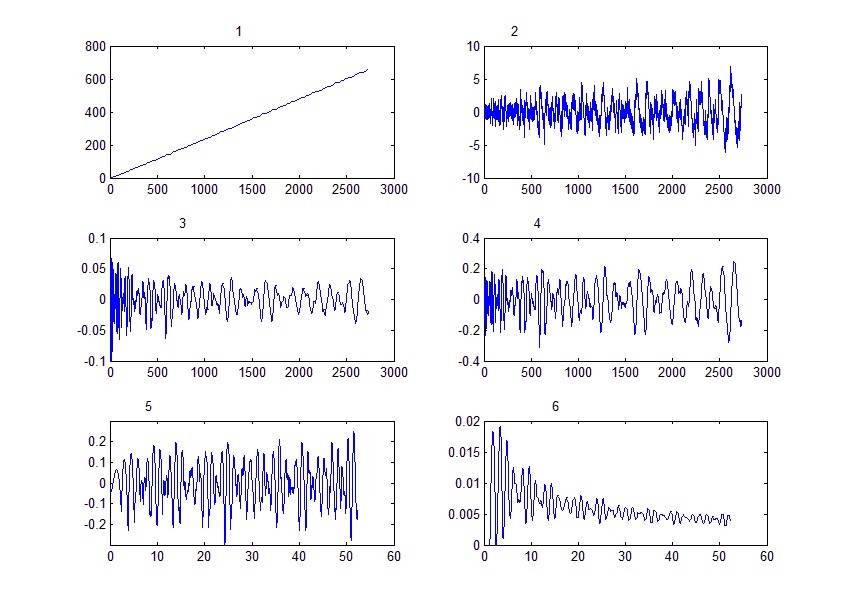}\\
\caption{Euclidean disk (Dirichlet conditions)}
\end{figure} 

%\begin{center}
%\includegraphics[scale=.4]{EucDDiscDifferenceGraphs}

%Figure 2.6
%\end{center} 

\subsection{Euclidean disk with Neumann boundary conditions (again with radius one)}

In this case $\lambda = (z_{n})^{2}$ where $J'_{k}(z_{n}) = 0$. Here we were able to obtain the exact values for the first 550 eigenvalues. Figure 4 shows the corresponding data with $\widetilde{N}(t) = \frac{1}{4}t + \frac{1}{2}\sqrt{t} + \frac{1}{6}$. The qualitative features observed for the previous example are evident here as well.

\begin{figure}\label{count_euc_n_disc}
\includegraphics[width=\scalor\textwidth]{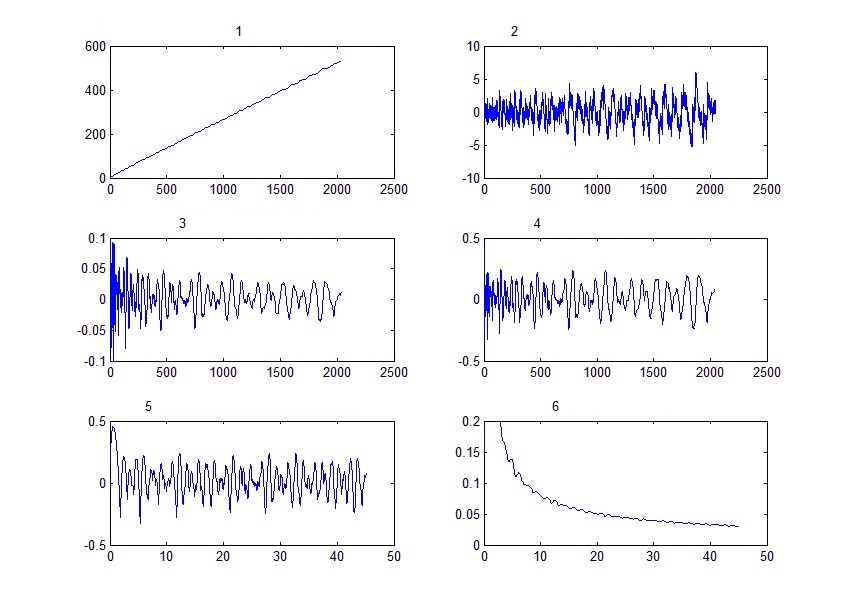}\\
\caption{Euclidean disk (Neumann conditions)}
\end{figure} 

%\begin{center}
%\includegraphics[scale=.6]{EucNDiscDifferenceGraphs}

%Figure 2.8
%\end{center} 

\section{Flat Surfaces}\label{sec_flat}
In this section we discuss examples of polygonal surfaces in Euclidean space. In particular we examined examples of nonconvex surfaces, surfaces with angles exceeding $\pi$, and surfaces that are not simply connected. There are still more examples on the website \cite{murray2015results}. For each example we give the counting function $\widetilde{N}(t)$ and the analog of Figure 1

\subsection{Triangle with Dirichlet boundary conditions}

The angles are $\theta_1 =\frac{\pi}{4}, \theta_2 = \frac{\pi}{5}, \theta_3 = \frac{11\pi}{20}$ and 
$$\widetilde{N}(t) = \frac{\sin{\theta_{2}}\sin{\theta_{1}}}{8\pi\sin{\theta_{3}}}t-\frac{\frac{\sin{\theta_1}}{\sin{\theta_3}}\sqrt{t}+\frac{\sin{\theta_2}}{\sin{\theta_3}}+1}{4\pi}+\frac{9}{22}$$
Here we used the first 130 calculated eigenvalues, as accuracy begins to break down after that point. The graphs in Figure 5 are analogous to those in Figure 1 and show similar behavior.

\begin{figure}[ht]\label{count_euc_tri2}
\includegraphics[width=\scalor\textwidth]{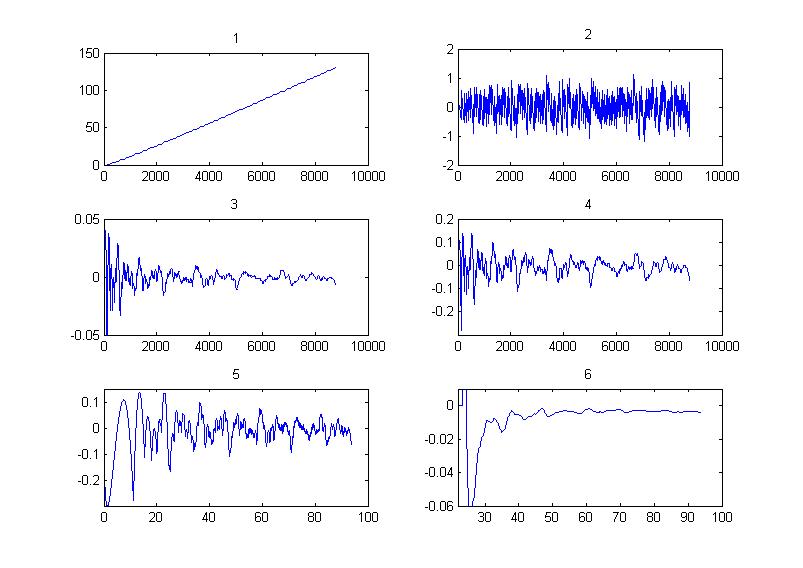}\\
\caption{Triangle with Dirichlet boundary conditions}
\end{figure} 

%\begin{center}
%\includegraphics[scale=.45]{EucDTri2DifferenceGraphs}
%Figure 3.2
%\end{center} 

\subsection{Triangle with Neumann boundary conditions}

This is the same triangle as above, with 
$$\widetilde{N}(t) = \frac{\sin{\theta_{2}}\sin{\theta_{1}}}{8\pi\sin{\theta_{3}}}+\frac{\frac{\sin{\theta_1}}{\sin{\theta_3}}t+\frac{\sin{\theta_2}}{\sin{\theta_3}}+1}{4\pi}\sqrt{t}+\frac{9}{22}$$
Here we used the first 150 calculated eigenvalues. The graphs in Figure 6 are analogous to Figure 5 and display the same behavior, except that in graph six of Figure 6 the graph is decreasing to 0, whereas it is increasing to zero in Figure 5. This difference is a result of the boundary conditions and is mirrored in all graphs of the same shape under Dirichlet and Neumann boundary conditions.
 
\begin{figure}[ht]\label{count_euc_n_tri2}
\centering
\includegraphics[width=\scalor\textwidth]{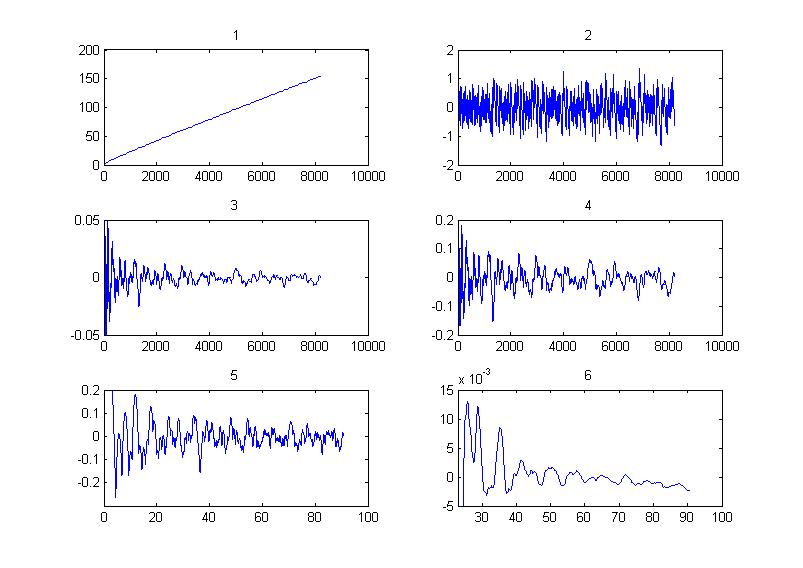}\\
\caption{Triangle with Neumann boundary conditions}
\end{figure} 

%\begin{center}
%\includegraphics[scale=.45]{EucNTri2DifferenceGraphs}\\
%Figure 3.4
%\end{center} 

\subsection{Triangle with mixed boundary conditions}

With the same triangle, we impose Dirichlet boundary conditions on $s_1$ and $s_2$ and Neumann boundary conditions on $s_3$. Here 

$$\widetilde{N}(t) = \frac{\sin{\theta_{2}}\sin{\theta_{1}}}{8\pi\sin{\theta_{3}}}t-\frac{\frac{\sin{\theta_1}}{\sin{\theta_3}}\sqrt{t}-\frac{\sin{\theta_2}}{\sin{\theta_3}}+1}{4\pi}\sqrt(t)+\frac{189}{110}$$
and the resulting graphs can be seen in Figure 7, which display similar behavior to those in Figures 5 and 6.

\begin{figure}[ht]\label{count_euc_m_tri2}
\centering
\includegraphics[width=\textwidth]{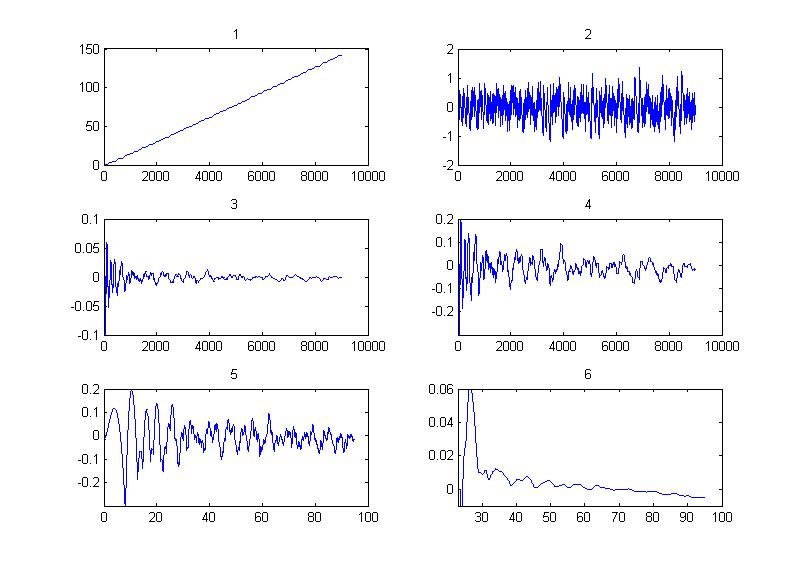}\\
\caption{Triangle with mixed boundary conditions}
\end{figure} 

%\begin{center}
%\includegraphics[scale=.45]{EucMixed2DifferenceGraphs}\\
%Figure 3.6
%\end{center} 

\subsection{Arrowhead with Dirichlet boundary conditions (see Figure 8)}

\begin{figure}[ht]\label{fig_arrowhead}
\centering
\includegraphics[width=\scalor\textwidth]{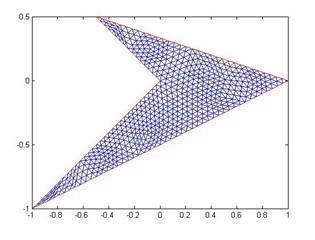}\\
\caption{Arrowhead region}
\end{figure} 

\iffalse We define an arrowhead to a Euclidean quadrilateral with one concave angle that is formed by joining two triangles with obtuse angles that have a pair of sides of equal length that are adjacent to the obtuse angle so that the concave angle of the resulting quadrilateral has a measure equal to the sum of the measures of the obtuse angles. The angle opposite the concave angle is the sum of the two non-obtuse angles adjacent to the sides that were joined together, and the remaining two angles are the remaining angles in each of the triangles. The sides of the arrowhead are then the sides of the triangles excluding the side they were joined on.\fi Here the sides are $s_1 = \frac{\sqrt{13}}{2}, s_2 = \frac{\sqrt{2}}{2}, s_3 = \sqrt{2}, s_4 = \sqrt{5}$. The angle $\theta_i$ joins sides $s_i$ and $s_{i-1}$; $\theta_1$ joins $s_1$ and $s_4$. The angle measures are $\theta_1 = \sin^{-1}{\frac{1}{\sqrt{13}}}+\sin^{-1}{\frac{1}{\sqrt{5}}}, \theta_2 = \cos^{-1}{(\frac{s^{2}_1+s^{2}_2-1}{2s_{1}s_{2}})},  \theta_4 = \cos^{-1}{(\frac{s^{2}_3+s^{2}_3-1}{2s_{3}s_{4}})}, \theta_3 = 2\pi-\theta_1-\theta_2-\theta_4$, and $\theta_3>\pi$. We then have 

$$a = \sqrt{\frac{s_1+s_2+1}{2}(s_1+\frac{s_2}{2}+1)(\frac{s_1}{2}+s_2+1)(s_1+s_2+\frac{1}{2})}$$

$$b = \sqrt{\frac{s_3+s_4+1}{2}(s_3+\frac{s_4}{2}+1)(\frac{s_3}{2}+s_4+1)(s_3+s_4+\frac{1}{2})}$$

and

$$\widetilde{N}(t) = \frac{a+b}{4\pi}t-\frac{\sum_{i=1}^{4}{s_i}}{4\pi}\sqrt{t}+\frac{\sum_{i=1}^{4}{\theta_i}}{24}$$

Figure 9 is analogous to the graphs we have seen before, and we can see that it displays the same behavior. Note that this region is not convex and contains an angle greater than $\pi$.

\begin{figure}[ht]\label{count_euc_arrowhead}
\includegraphics[width=\scalor\textwidth]{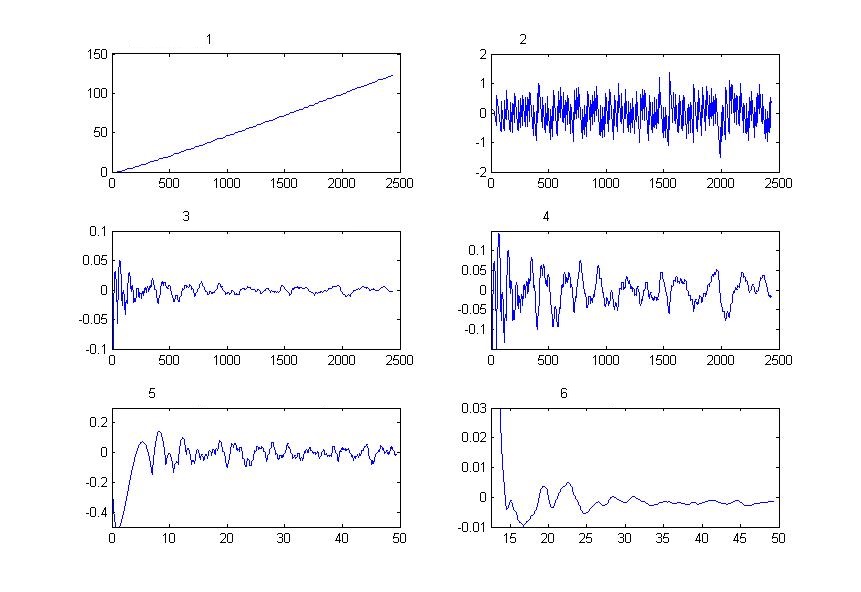}\\
\caption{Arrowhead with Dirichlet boundary conditions}
\end{figure} 

%\begin{center}
%\includegraphics[scale=.40]{mEucDArr1DifferenceGraphs}\\
%Figure 3.9
%\end{center} 

\subsection{Region between triangles with Dirichlet boundary conditions (see Figure 10)}

\begin{figure}[ht]\label{fig_euc_tri_in_tri}
\centering
\includegraphics[width=\scalor\textwidth]{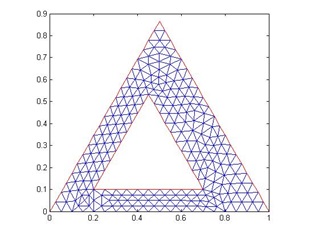}\\
\caption{Region between triangles}
\end{figure} 

This is not simply connected. The angles of the interior triangle are viewed from the surface and hence are the exterior angles. We need to keep the vertices of the inner triangle a reasonable distance from the edges of the outer triangle in order to have reasonable accuracy in computing eigenvalues. Note that the formula $\widetilde{N}(t) = \frac{3\sqrt{3}}{64\pi}t-\frac{9}{8\pi}{\sqrt{t}}+\frac{1}{15}$ is independent of the location and orientation of the inner triangle. The graphs seen in Figure 11 are analagous to those seen before and display similar results.

\begin{figure}[ht]\label{count_euc_tri_in_tri}
\centering
\includegraphics[width=\scalor\textwidth]{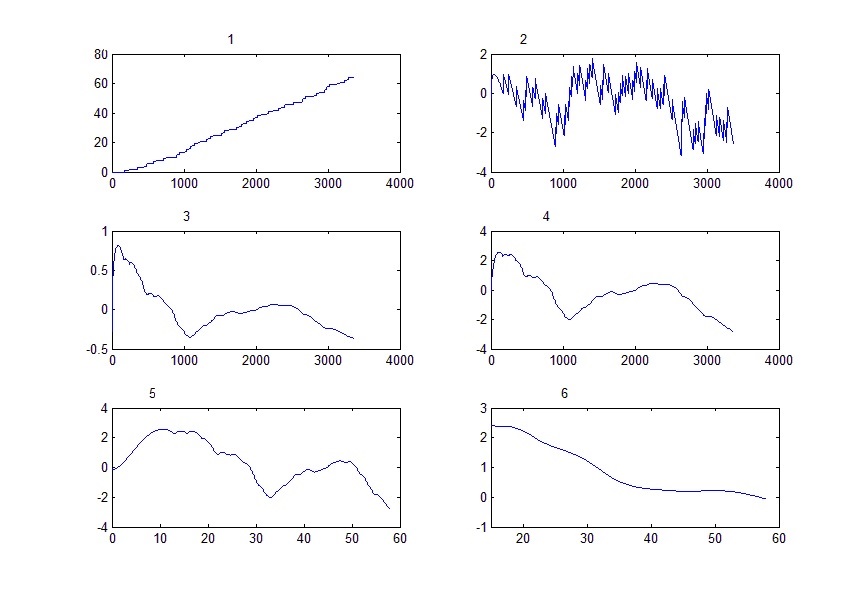}\\
\caption{Region between triangles with Dirichlet boundary conditions}
\end{figure} 

%\begin{center}
%\includegraphics[scale=.40]{EucDTrinTri1DifferenceGraphs}\\
%Figure 3.12
%\end{center} 

\subsection{Regular pentagon with Dirichlet boundary conditions}

Here $\widetilde{N}(t) = \frac{\sqrt{5(5+2\sqrt{5}}}{16\pi}t-\frac{5}{4\pi}\sqrt{t}+\frac{2}{9}$. Note that there are many eigenvalues of multiplicity two, due to the $D_5$ symmetry group. This gives us a reasonable tool for assessing the accuracy of our computations (since MATLAB does not select symmetric triangulations). The first 10 eigenvalues are displayed in Table 3.

\begin{table}[ht]\label{tab_euc_pent}
\centering
	\begin{tabular}{| l | l | l | l | l | l | l | l | }
	\hline
	 & Initial & 1 & 2 & 3 & 4 & 5 & Predicted \\ \hline
1 & 11.1479 & 11.03526 & 11.00624 & 10.99889 & 10.99704 & 10.99658 & 10.99643 \\ \hline
2 & 28.8006 & 28.04273 & 27.85074 & 27.80238 & 27.79025 & 27.78721 & 27.7862 \\ \hline
3 & 28.8171 & 28.04763 & 27.85201 & 27.8027 & 27.79033 & 27.78723 & 27.7862 \\ \hline
4 & 52.4110 & 50.06061 & 49.47096 & 49.323 & 49.28594 & 49.27667 & 49.27358 \\ \hline
5 & 52.5710 & 50.10611 & 49.48295 & 49.32606 & 49.28671 & 49.27686 & 49.27359 \\ \hline
6 & 61.5873 & 58.21567 & 57.37533 & 57.16477 & 57.11204 & 57.09885 & 57.09447 \\ \hline
7 & 84.9386 & 78.95964 & 77.47257 & 77.1007 & 77.00764 & 76.98437 & 76.97664 \\ \hline
8 & 85.0667 & 78.9798 & 77.47749 & 77.10193 & 77.00795 & 76.98445 & 76.97664 \\ \hline
9 & 99.4387 & 91.724 & 89.80689 & 89.32718 & 89.2071 & 89.17706 & 89.16708 \\ \hline
10 & 100.119 & 91.89246 & 89.84908 & 89.33775 & 89.20975 & 89.17773 & 89.16708 \\ \hline
	\end{tabular}

	\caption{Regular pentagon with Dirichlet boundary conditions}
\end{table}

\begin{figure}[ht]\label{count_euc_pent}
\includegraphics[width=\scalor\textwidth]{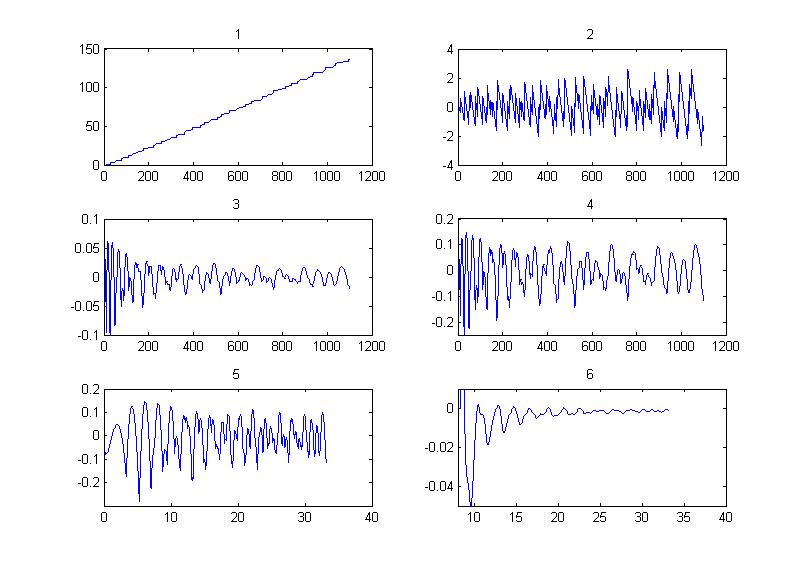}\\
\caption{Regular pentagon with Dirichlet boundary conditions}
\end{figure} 

%\begin{center}
%\includegraphics[scale=.37]{EucDPentDifferenceGraphs}\\
%Figure 3.14
%\end{center} 

\subsection{Regular hexagon with Dirichlet boundary conditions}

Here $\widetilde{N}(t) = \frac{3\sqrt{3}}{8\pi}t-\frac{3}{\pi}\sqrt{t}+\frac{5}{24}$. Note that all Dirichlet eigenfunctions of the equilateral triangle extend by odd reflections to Dirichlet eigenfunctions of the hexagon with the same eigenvalue. In our table of eigenvalues we therefore divide by $(\frac{4}{3}\pi)^{2}$ so that these eigenvalues become integers. This gives us an accuracy check. We have a $D_{6}$ symmetry group so that most eigenvalues have multiplicity two.

\begin{table}[ht]\label{tab_euc_hex}
\centering
	\begin{tabular}{| l | l | l | l | l | l | l | l | }
	\hline
	 & Initial & 1 & 2 & 3 & 4 & 5 & Predicted \\ \hline
1 & 0.413621 & 0.409345 & 0.408200 & 0.407905 & 0.40783 & 0.40781 & 0.40781 \\ \hline
2 & 1.069612 & 1.042747 & 1.035759 & 1.03398 & 1.03353 & 1.03342 & 1.03338 \\ \hline
3 & 1.07095 & 1.043083 & 1.035843 & 1.034001 & 1.03354 & 1.03342 & 1.03338 \\ \hline
4 & 1.968578 & 1.880006 & 1.857227 & 1.851466 & 1.85002 & 1.84966 & 1.84953 \\ \hline
5 & 1.969445 & 1.880229 & 1.857298 & 1.851485 & 1.85002 & 1.84966 & 1.84953 \\ \hline
6 & 2.288059 & 2.175413 & 2.146518 & 2.139204 & 2.13737 & 2.13691 & 2.13675 \\ \hline
7 & 2.976487 & 2.782209 & 2.731673 & 2.718854 & 2.71563 & 2.71482 & 2.71455 \\ \hline
8 & 3.328669 & 3.08177 & 3.020466 & 3.005122 & 3.00128 & 3.0003 & 3 \\ \hline
9 & 3.808021 & 3.523318 & 3.450285 & 3.431787 & 3.42714 & 3.4260 & 3.42558 \\ \hline
10 & 3.816109 & 3.524926 & 3.450607 & 3.43186 & 3.42715 & 3.4260 & 3.42558 \\ \hline
	\end{tabular}

	\caption{Regular hexagon with Dirichlet boundary conditions}
\end{table}

\begin{figure}[ht]\label{count_euc_hex}
\centering
\includegraphics[width=\scalor\textwidth]{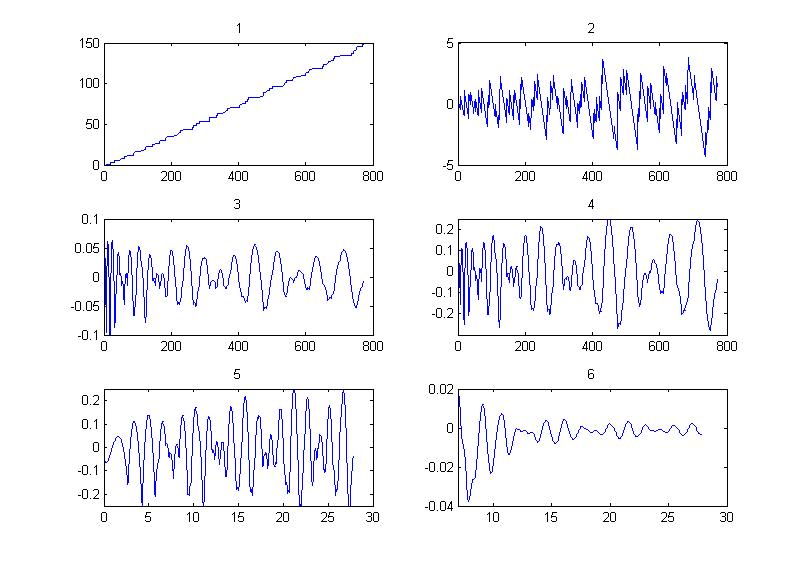}\\
\caption{Regular hexagon with Dirichlet boundary conditions}
\end{figure} 

%\begin{center}
%\includegraphics[scale=.40]{EucDHexDifferenceGraphs}\\
%Figure 3.16
%\end{center} 

\subsection{6-regular star with Dirichlet boundary conditions}

Here $\widetilde{N}(t) = \frac{s\sqrt{3}}{4\pi}t -\frac{3}{\pi}\sqrt{t}+\frac{25}{48}$. As in the case of the hexagon, Dirichlet eigenfunctions of the equilateral triangle extend by odd reflection, and there is a $D_6$ symmetry group. Therefore we again divide our table by $(\frac{4}{3}\pi)^{2}$ so that these eigenvalues become integers. The other eigenvalues of the hexagon do not, however, extend to the 6-regular star.

\begin{figure}[ht]\label{count_euc_star}
\centering
\includegraphics[width=\scalor\textwidth]{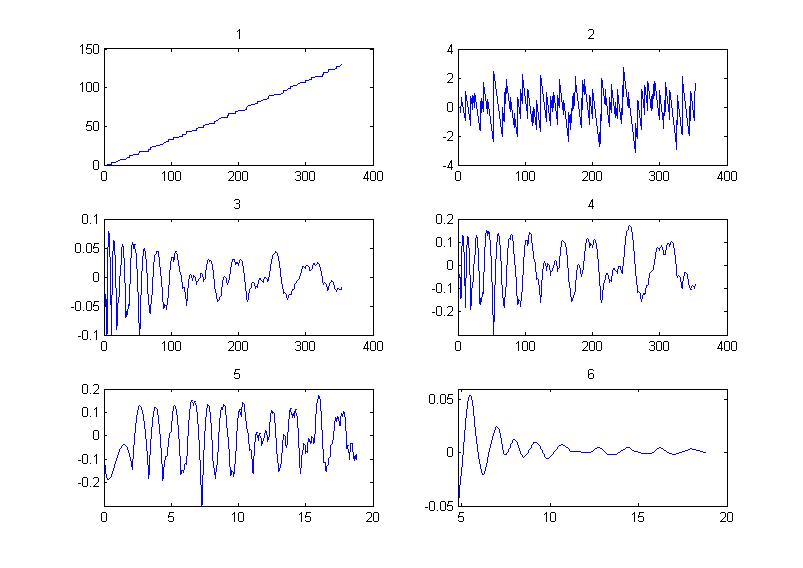}\\
\caption{6-regular star with Dirichlet boundary conditions}
\end{figure} 

%\begin{center}
%\includegraphics[scale=.40]{EucDSixStarDifferenceGraphs}\\
%Figure 3.18
%\end{center} 

\section{Hyperbolic Surfaces}\label{sec_hyp}

In this section we discuss examples of surfaces in the hyperbolic plane of constant negative curvature $-1$. We use the upper half-plane model. In this model the Laplacian is given by 

\begin{equation}
\tag{4.1} \Delta = y^2(\frac{\partial^2}{\partial x^2}+\frac{\partial^2}{\partial y^2})
\end{equation}

so the eigenvalue problem

\begin{equation}
\tag{4.2} -\Delta u = \lambda u
\end{equation}

is transformed into

\begin{equation}
\tag{4.3} -(\frac{\partial^2}{\partial x^2}+\frac{\partial^2}{\partial y^2})u(x,y) = \lambda y^{2}u(x,y)
\end{equation}

and we used MATLAB to solve $(4.3)$ on the surfaces with the appropriate boundary conditions. For simplicity we restricted our attention to Dirichlet boundary conditions, and our surfaces were either disks or triangles.

To describe triangles we recall that geodesics in the upper half-plane model are either vertical half lines or half circles that intersect the x-axis perpendicularly. Without loss of generality we may take one side of the triangle to lie along the y-axis. Specifically, the triangle will have vertices $(0,y_1), (0,y_2)$ and $(x_3,y_3)$, seen in Figure 15 as points C, A, and B, respectively.
\begin{figure}[ht]\label{fig_hyp_tri}
\centering
\includegraphics[width=0.5\textwidth]{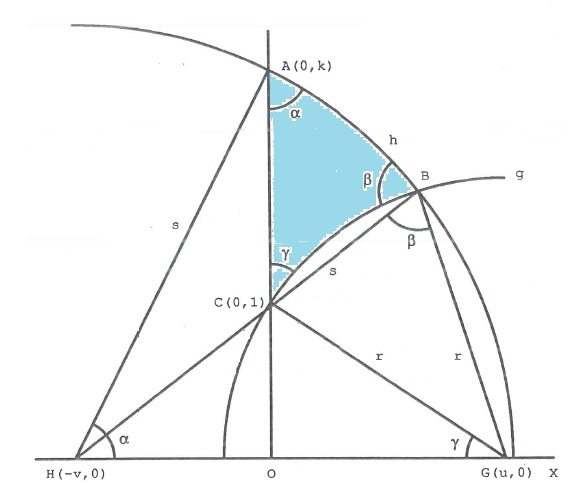}\\
\caption{Hyperbolic Triangle}
\end{figure} 

%here

The two boundary circles are $y^2 + (x-a_j)^2 = r_{j}^2$ for $j = 1,2$, and $(x_3,y_3)$ lies at the intersection of these circles, so $x_3 = \frac{r_2^2-r_1^2-a_2^2+a_1^2}{s(a_1-a_2)}, y_3 = \sqrt{r_1^2-(x_3-a_1)^2}$, and also $y_j = \sqrt{r_j^2-a_j^2}$ for $j = 1,2$.

Since the model is conformal, the angles are the same as the Euclidean angles, so we have $\alpha_1 = \frac{\pi}{2}-\tan^{-1}{(\frac{a_1}{y_1})}, \alpha_2 = \frac{\pi}{2} - \tan^{-1}{(\frac{-a_2}{y_2})}$ and \newline{} $\alpha_3 = \tan^{-1}{(\frac{a_1-x_3}{y_3})}+\tan^{-1}{(\frac{x_3-a_a}{y_3})}$. The lengths of the opposite sides are $L_j = \frac{1}{2}log(\frac{r_j+a_j}{r_j-a_j})-\frac{1}{2}log(\frac{r_j-x_x+a_j}{r_j+x_3+a_j})$ for $j = 1,2$ and $L_3 = log\frac{y_2}{y_1}$. The area of the triangle is 
$$A=\int\!\!\!\int_T \frac{dxdy}{y^2}=\frac{1}{r_2}(\cos^{-1}(\frac{-a_2}{r_2})-\cos^{-1}(\frac{x_3-a_2}{r_2}))+\frac{1}{r_1}(\cos^{-1}(\frac{x_3-a_1}{r_1})-\cos^{-1}(\frac{-a_1}{r_1}))$$ Thus we have

\begin{equation}
\tag{4.4} \widetilde{N}(t)=\frac{1}{4\pi}At-\frac{1}{4\pi}(L_1+L_2+L_3)t^{\frac{1}{2}}+C
\end{equation}

for

\begin{equation}
\tag{4.5} C=-\frac{1}{12\pi}A+\frac{1}{24}\sum_{j=0}^{2}(\frac{\pi}{\alpha_j}-\frac{\alpha_j}{\pi})
\end{equation}

Of course everything may be expressed entirely in terms of the angles, since the angles determine the triangle. Thus the hyperbolic law of cosines yields

\begin{equation}
\tag{4.6} L_i=\cosh^{-1}(\frac{\cos\alpha_j\cos\alpha_k+\cos\alpha_i}{\sin\alpha_j\sin\alpha_k})
\end{equation}

for $(i,j,k)$ any permutation of $(1,2,3)$, and the angle defect formula yields

\begin{equation}
\tag{4.7} A=\pi-(\alpha_1+\alpha_2+\alpha_3)
\end{equation}

\subsection{Hyperbolic Equilateral Triangles with Dirichlet boundary conditions}

We take $\alpha_1 = \alpha_2 = \alpha_3 = \frac{\pi}{k}$ for $k$ an integer, $k \geq 4$. These triangles tesselate the hyperbolic plane. When $k$ is even we may take odd reflections of the Dirichlet eigenfunctions to see that we are generating a subset of the collection of eigenfunctions on the hyperbolic closed manifolds $\Gamma \backslash SL(2,\mathbb{R})/SO(2)$ for the appropriate discrete subgroup $\Gamma$.

We show the results for $k = 4,6$ in Figures 16 and 17, respectively. Already for $k = 6$ the accuracy of our approximations begins to degrade. The website shows complete data for $k = 4,5,6,7$.

\begin{figure}[ht]\label{count_hyp_equ_tri4}
\centering
\includegraphics[width=\scalor\textwidth]{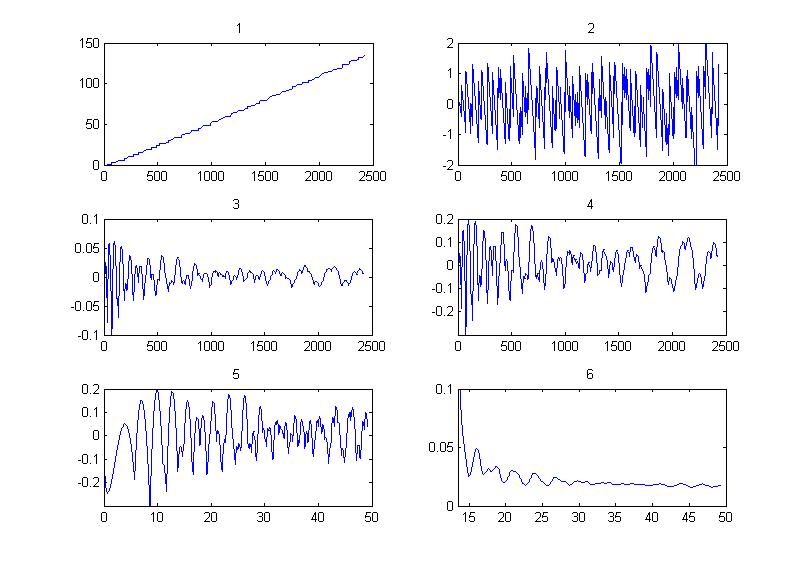}\\
\caption{Hyperbolic equilateral triangle with $k = 4$}
\end{figure} 

%\begin{center}
%\includegraphics[scale=.4]{HypDEquTri4DifferenceGraphs}\\
%Figure 4.3
%\end{center} 

\begin{figure}[ht]\label{count_hyp_equ_tri6}
\centering	
\includegraphics[width=\textwidth]{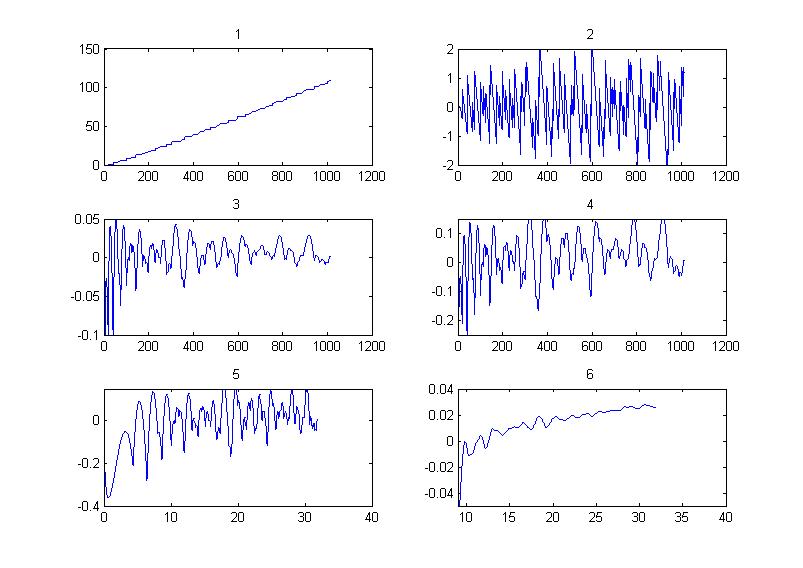}\\
\caption{Hyperbolic equilateral triangle with $k = 6$}
\end{figure} 

%\begin{center}
%\includegraphics[scale=.40]{HypDEquTri6DifferenceGraphs}\\
%Figure 4.5
%\end{center} 

\subsection{General Hyperbolic Triangles with Dirichlet boundary conditions}

We present two hyperbolic triangles here with arbitrary measurements. Triangles are specified by a label $(u,v,s)$ which correspond to the measurements in Figure 15. Figure 18 corresponds to $(5,10,11)$ and Figures 19 corresponds to $(3,4,6)$. Note that in the fifth and sixth counting graphs in figures 18 and 19 we begin to lose accuracy more quickly than we do in the Euclidean results. This is not unique to the arbitrary triangles, as it is present in both the hyperbolic equilateral triangles and hyperbolic discs, but it is especially noticeable here.

The complete results for more arbitrary hyperbolic triangles are shown on the website.

\begin{figure}[ht]\label{count_hyp_tri1}
\centering
\includegraphics[width=\scalor\textwidth]{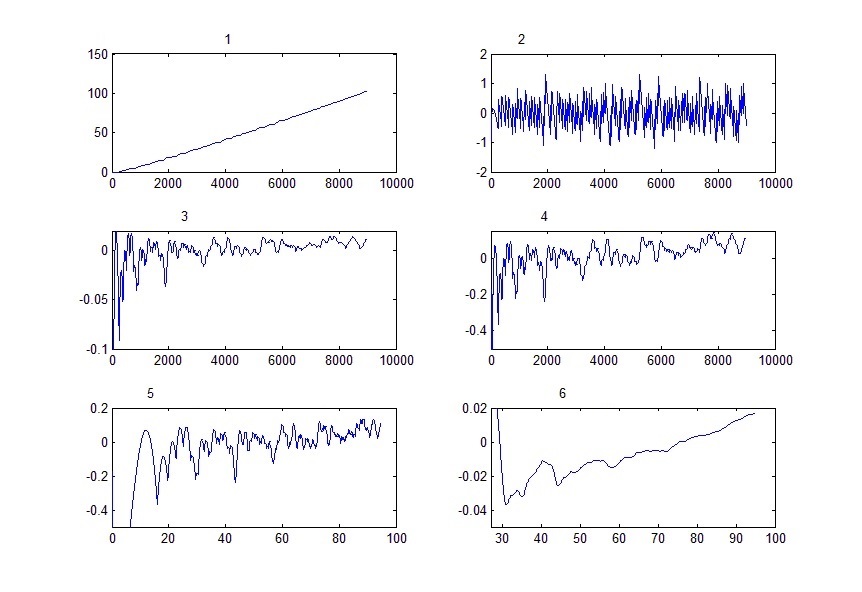}\\
\caption{First hyperbolic triangle}
\end{figure} 

%\begin{center}
%\includegraphics[scale=.5]{HypDTri1DifferenceGraphs}\\
%Figure 4.7
%\end{center} 

\begin{figure}[ht]\label{count_hyp_tri2}
\centering
\includegraphics[width=\scalor\textwidth]{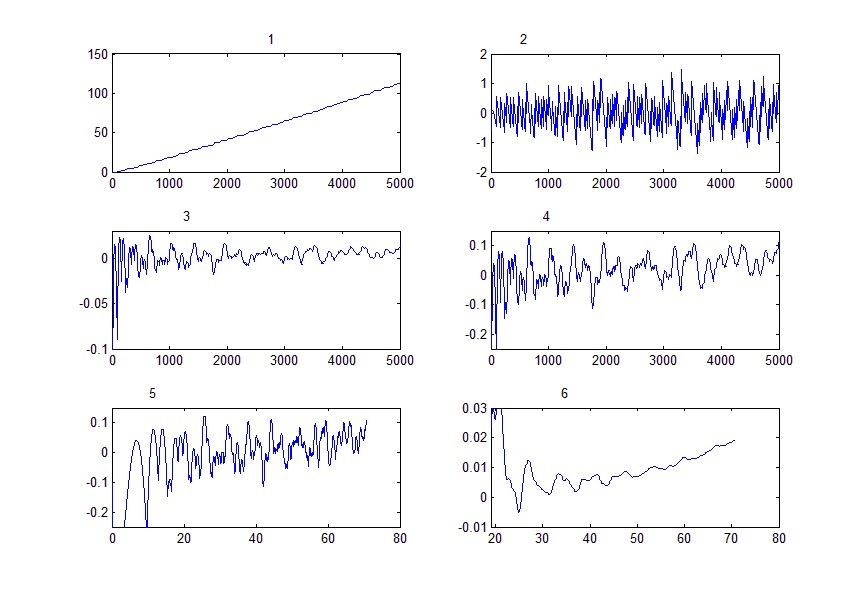}\\
\caption{Second hyperbolic triangle}
\end{figure} 

%\begin{center}
%\includegraphics[scale=.5]{HypDTri2DifferenceGraphs}\\
%Figure 4.9
%\end{center} 

\subsection{Hyperbolic Discs with Dirichlet boundary conditions}

We are able to calculate eigenvalues on a disc of hyperbolic radius $ R$ by calculating eigenvalues for a Euclidean disc of radius $r=\frac{e^{2R}}{2}$ centered at $\frac{e^{2R}}{2}+1$. The resulting disc has area $A=4\pi\sinh^{2}(\frac{R}{2})$ and circumference $C = 2\pi\sinh(R)$ We then have $\widetilde{N}(t)=\frac{A}{4\pi}t-\frac{C}{r\pi}\sqrt{t}+\frac{1}{6}$. We can see from figures 20 and 21 that this $\widetilde{N}(t)$ appears to be strongly supported, though as is the case of the hyperbolic triangles, we begin to lose accuracy in the predicted eigenvalues more quickly here than in the Euclidean case. This becomes a particular issue for the MATLAB PDE solver in the case of discs however, as the radius of the Euclidean disc we solve for eigenvalues on grows exponentially with the hyperbolic radius, leading to longer computation times as mesh with an exponentially growing number of points is needed to estimate the values. The discs in Figures 20 and 21 have R = 1, 1/2, respectively.

\begin{figure}[ht]\label{count_hyp_disc1}
\centering
\includegraphics[width=\scalor\textwidth]{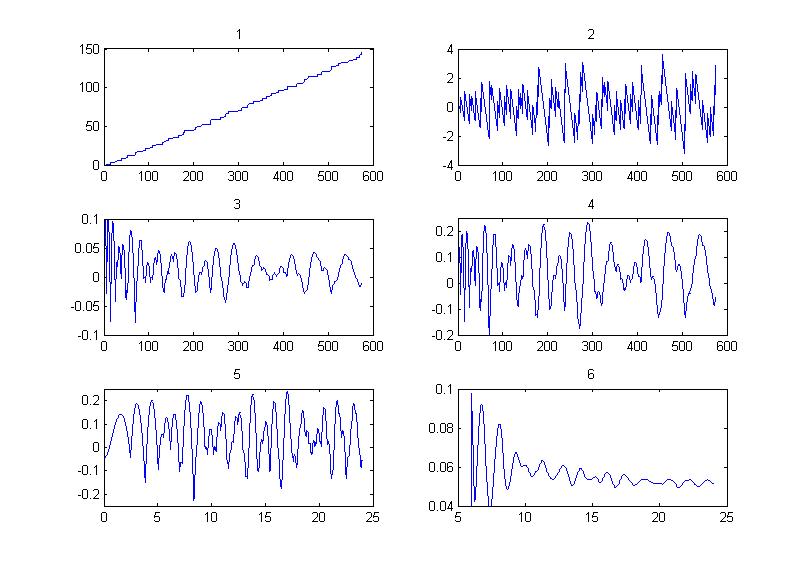}\\
\caption{Hyperbolic disk with $R = 1$}
\end{figure} 

%\begin{center}
%\includegraphics[scale=.40]{HypDDisc1DifferenceGraphs}\\
%Figure 4.11
%\end{center} 

\begin{figure}[ht]\label{count_hyp_disc_half}
\centering
\includegraphics[width=\scalor\textwidth]{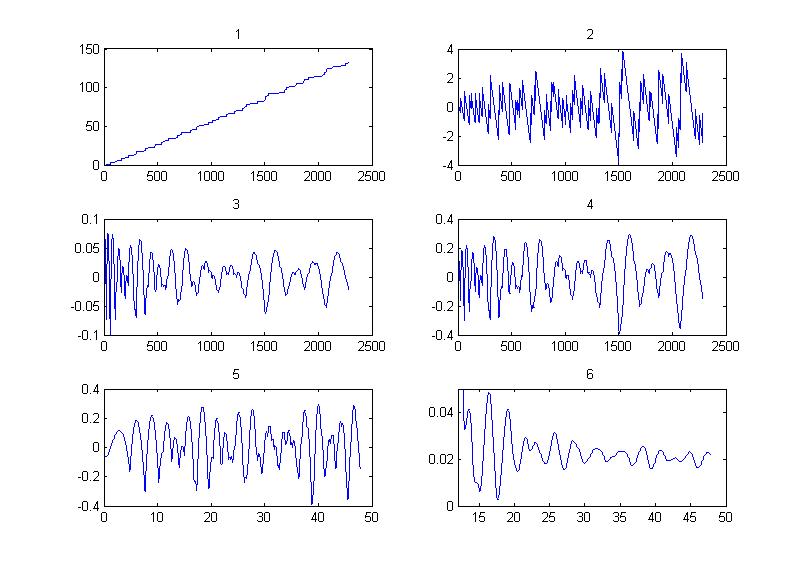}\\
\caption{Hyperbolic disk with $R = \frac{1}{2}$}
\end{figure} 

%\begin{center}
%\includegraphics[scale=.40]{HypDDiscHalfDifferenceGraphs}\\
%Figure 4.13
%\end{center} 

\section{Spherical Surfaces}\label{sec_sph}

In this section we discuss examples of surfaces in the unit sphere (curvature +1). We use stereographic projection, placing the center of the sphere at $(0,0,1)$ and projecting from $(0,0,2)$ onto the $(u,v)$ plane by $u = \frac{2x}{2-z}, v = \frac{2y}{2-z}$. The equator is mapped to the circle $u^2+v^2=4$, great circles through the poles are mapped to the lines through the origin, and other great circles are mapped to circles intersecting $u^2+v^2=4$ at two antipodal points. If we fix parameters to $t>0$ and $\beta$ then these circles are given by $(u-t\sin{\beta})^2+(v+t\cos{\beta})^2=t^2+4$ (intersecting $u^2+v^2=4$ at $\pm (2\cos{\beta},2\sin{\beta})$).

We will consider triangles with vertices $(u_1,0), (u_2,0)$ and $(u_3,v_3)$, with one edge along the $u$-axis and two edges being arcs of circles $(u-t_j\sin{\beta_j})^2+(v+t_j\cos{\beta_j})^2=t_j^2+4$ for $j=1,2$. The angles of the triangle are given by
$$\alpha_1={\tan}^{-1}{(\frac{-u_1+t_1\sin{\beta_1}}{t_1\cos{\beta_1}})}$$
$$\alpha_2=\pi-{\tan}^{-1}{(\frac{-u_2+t_2\sin{\beta_2}}{t_2\cos{\beta_2}})}$$
$$\alpha_3={\tan}^{-1}{(\frac{-u_3+t_2\sin{\beta_2}}{v_3+t_1\cos{\beta_2}})-{\tan}^{-1}{(\frac{-u_3+t_1\sin{\beta_1}}{v_3+t_1\cos{\beta_1}})}}$$
The angles completely determine the triangle. The lengths of the sides are given by the spherical law of cosines\\
\begin{equation}
\tag{5.1} L_i= \cos^{-1}{(\frac{\cos{\alpha_i}+\cos{\alpha_j}\cos{\alpha_k}}{\sin{\alpha_j}\sin{\alpha_k}})}\\
\end{equation}
for $(i,j,k)$ a permutation of $(1,2,3)$, and the area is given by the angle defect\\
\begin{equation}
\tag{5.2} A = (\alpha_1+\alpha_2+\alpha_3) - \pi
\end{equation}
The Laplacian is given by
\begin{equation}
\tag{5.3} \Delta = (\frac{u^2+v^2+4}{4})^2(\frac{\partial^2}{\partial u^2}+\frac{\partial^2}{\partial v^2})
\end{equation}

\subsection{Spherical Equilateral Right Triangle with Dirichlet boundary conditions}

This triangle serves as our main accuracy check for our calculated eigenvalues in spherical space. This is because it is a region for which the eigenvalue spectrum is known: the $i^{th}$ distinct eigenvalue is equal to $4i^2+6i+2$ and has multiplicity $i$. We can see the first eigenvalues in Table 5 below.

\begin{table}[ht]\label{tab_sph_equ_right_tri}
	\centering
	\begin{tabular}{| l | l | l | l | l | l | l | l | l | }
	\hline
	 & Initial & 1 & 2 & 3 & 4 & 5 & Predict & True\\ \hline
1 & 12.1683 & 12.0426 & 12.0107 & 12.0027 & 12.0007 & 12.0002 & 12 & 12 \\ \hline
2 & 30.9285 & 30.2355 & 30.0593 & 30.0148 & 30.0037 & 30.0009 & 30 & 30 \\ \hline
3 & 31.1082 & 30.2803 & 30.0704 & 30.0176 & 30.0044 & 30.0011 & 30 & 30 \\ \hline
4 & 58.8956 & 56.7339 & 56.1845 & 56.0462 & 56.0116 & 56.0029 & 56 & 56 \\ \hline
5 & 59.6055 & 56.9108 & 56.2287 & 56.0573 & 56.0143 & 56.0036 & 56.0000 & 56 \\ \hline
6 & 59.8717 & 56.9775 & 56.2454 & 56.0615 & 56.0154 & 56.0038 & 56.0000 & 56 \\ \hline
7 & 97.1655 & 91.8111 & 90.4552 & 90.1140 & 90.0285 & 90.0071 & 90.0000 & 90 \\ \hline
8 & 98.8889 & 92.2347 & 90.5603 & 90.1402 & 90.0351 & 90.0088 & 90.0000 & 90 \\ \hline
9 & 99.7452 & 92.4530 & 90.6150 & 90.1539 & 90.0385 & 90.0097 & 90.0000 & 90 \\ \hline
10 & 100.1208 & 92.5527 & 90.6404 & 90.1603 & 90.0401 & 90.0100 & 90.0000 & 90 \\ \hline
	\end{tabular}

	\caption{Spherical equilateral right triangle}
\end{table}

In the graphical data displayed in Figure 22 (and subsequent figures in this section) we show
\begin{enumerate}
\item $N(t)$
\item $D(t) = N(t) - \widetilde{N}(t)$
\item $A(t) = \frac{1}{t}\int_{0}^{t}D(s)ds$
\item $A(t^{2})$
\item $\frac{1}{t-a}\int_{a}^{t}s^{\frac{1}{2}}A(s^{2})ds$ for $a = $ the highest predicted eigenvalue divided by 16, removing the first $\frac{1}{4}$ of graph 4 from figuring into graph 5 and eliminating potential early extreme values so that it converges to 0 more quickly.
\end{enumerate}
The scales for each of these graphical displays are the same as the scales in the analagous set of six graphs we used for Euclidean and hyperbolic regions.

\begin{figure}[ht]\label{count_sph_equ_right_tri}
\centering
\includegraphics[width=\scalor\textwidth]{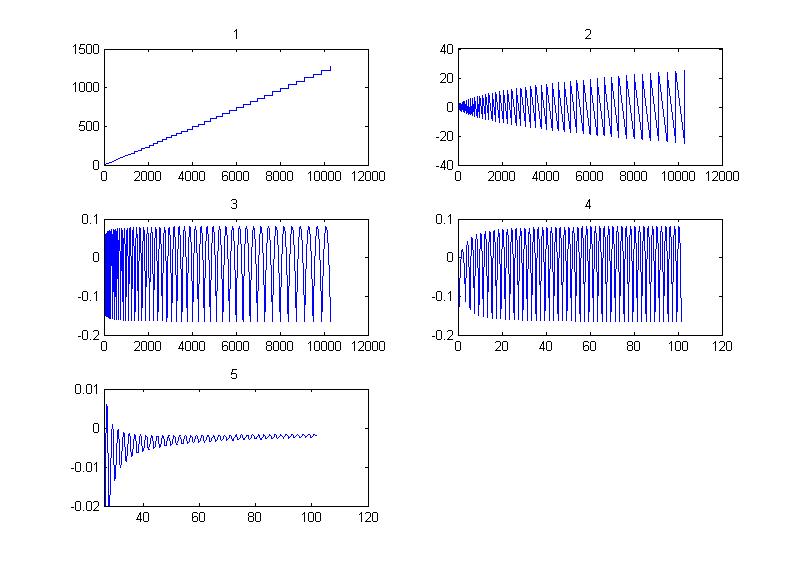}\\
\caption{Spherical equilateral right triangle}
\end{figure} 

%\begin{center}
%\includegraphics[scale=.44]{SphEquRightDTriDifferenceGraphs}\\
%Figure 5.2
%\end{center} 

\subsection{General Spherical Triangle with Dirichlet boundary conditions}

We will now present a spherical triangle with arbitrary measurements. Spherical triangles are specified by a label $(t_1,\beta_1,t_2,\beta_2)$ which correspond to the measurements described above. Here we have $\widetilde{N}(t)=\frac{Area}{4\pi}t-\frac{\Sigma_{i=1}^3{L_i}}{4\pi}\sqrt{t}+(\frac{Area}{12\pi}+\frac{1}{24}\Sigma_{i=1}^3{(\frac{\pi}{\alpha_i}-\frac{\alpha_i}{\pi})})$ The triangle corresponding to Figure 23 is $(-1.5,\frac{\pi}{4},-2,-\frac{\pi}{6})$. Note in the figure that while the accuracy of our calculated eigenvalues suffers some decay, it does so at a slower rate than in the Hyperbolic surfaces.

\begin{figure}[ht]\label{count_sph_tri2}
\centering
\includegraphics[width=\scalor\textwidth]{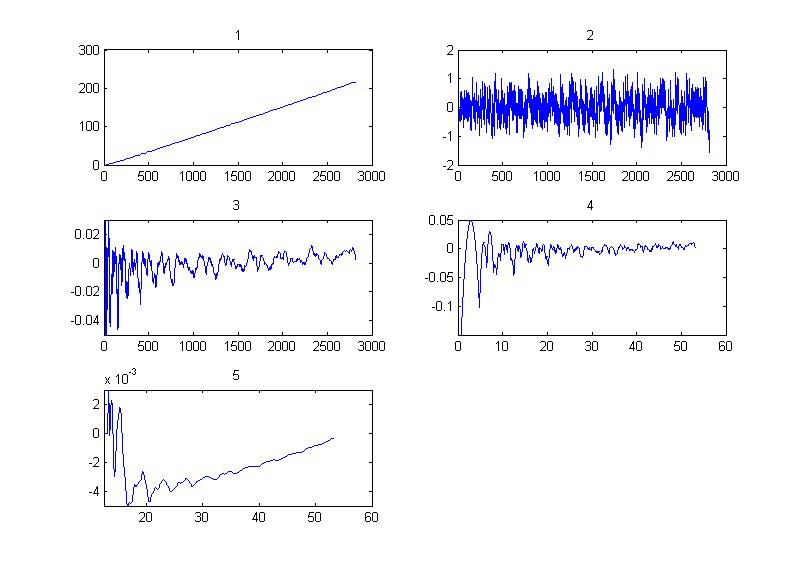}\\
\caption{General spherical triangle}
\end{figure} 

%\begin{center}
%\includegraphics[scale=.44]{SphDTri2DifferenceGraphs}\\
%Figure 5.4
%\end{center} 

\subsection{Spherical Disc with Dirichlet boundary conditions}

The final case we wish to present is that of the spherical disc. As with the spherical triangles, we use stereographic projection to solve for eigenvalues on a Euclidean disc using (5.3) as the Laplacian. A Euclidean disc centered at the origin with radius $r$ corresponds to a spherical disc with radius $R=\frac{2\sin{r}}{1+\cos{r}}$. The spherical disc then has area $A=2\pi(1-\cos{r})$ and circumference $C=2\pi\sin{r}$. This gives us $\widetilde{N}(t)=\frac{A}{4\pi}t-\frac{C}{4\pi}\sqrt{t}+\frac{1}{6}$. Here the hemisphere, $r=\frac{\pi}{2}$, serves as a test case for the accuracy of our predicted eigenvalues: it has a known eigenvalue spectrum such that the $n^{th}$ unique eigenvalue is $n(n+1)$ and has multiplicity $n$. Table 6 and Figures 24 strongly support the accuracy of $\widetilde{N}(t)$. Figure 25 is the corresponding graphs for the disc with $r=\frac{\pi}{4}$, and likewise strongly supports  our predicted eigenvalues. Unfortunately, as can be seen in graphs 3, 4, and 5 of 25, our accuracy again begins to drop after a point, but they are accurate enough to support our calculation of $\widetilde{N}(t)$.

\begin{table}[ht]\label{tab_sph_hemisphere}
	\centering
	\begin{tabular}{| l | l | l | l | l | l | l | l | l | }
	\hline
	 & Initial & 1 & 2 & 3 & 4 & 5 & Predict & True\\ \hline
1 & 2.03144 & 2.00787 & 2.00197 & 2.00049 & 2.00012 & 2.00003 & 2 & 2 \\ \hline
2 & 6.19516 & 6.04907 & 6.0123 & 6.00308 & 6.00077 & 6.00019 & 6 & 6 \\ \hline
3 & 6.20475 & 6.05152 & 6.0129 & 6.00324 & 6.00080 & 6.00020 & 6 & 6 \\ \hline
4 & 12.7272 & 12.1822 & 12.0457 & 12.0114 & 12.0029 & 12.0007 & 12.0000 & 12 \\ \hline
5 & 12.7330 & 12.1830 & 12.0458 & 12.0115 & 12.0029 & 12.0007 & 12.0000 & 12 \\ \hline
6 & 12.8645 & 12.2143 & 12.0535 & 12.0134 & 12.0034 & 12.0008 & 12.0000 & 12 \\ \hline
7 & 21.8055 & 20.4515 & 20.1130 & 20.0283 & 20.0071 & 20.0018 & 20.0000 & 20 \\ \hline
8 & 21.8512 & 20.4617 & 20.1155 & 20.0289 & 20.0072 & 20.0018 & 20.0000 & 20 \\ \hline
9 & 22.2046 & 20.5453 & 20.1362 & 20.0341 & 20.0085 & 20.0021 & 20.0000 & 20 \\ \hline
10 & 22.3183 & 20.5798 & 20.1451 & 20.0363 & 20.0091 & 20.0023 & 20.0000 & 20 \\ \hline
252 & 0 & 0 & 0 & 530.180 & 512.023 & 507.504 & 505.985 & 506 \\ \hline
253 & 0 & 0 & 0 & 530.425 & 512.092 & 507.522 & 505.986 & 506 \\ \hline
	\end{tabular}

	\caption{Hemisphere}
\end{table}

\begin{figure}[ht]\label{count_sph_hemisphere}
\centering
\includegraphics[width=\scalor\textwidth]{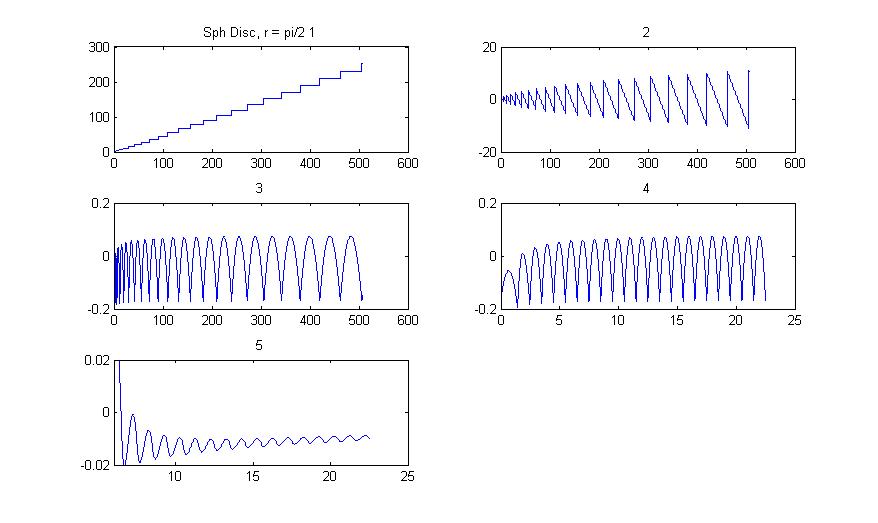}\\
\caption{Hemisphere}
\end{figure} 

%\begin{center}
%\includegraphics[scale=.40]{SphDisc_0_0DifferenceGraphs}\\
%Figure 5.6
%\end{center} 

\begin{figure}[ht]\label{count_sph_disc2}
\centering
\includegraphics[width=\scalor\textwidth]{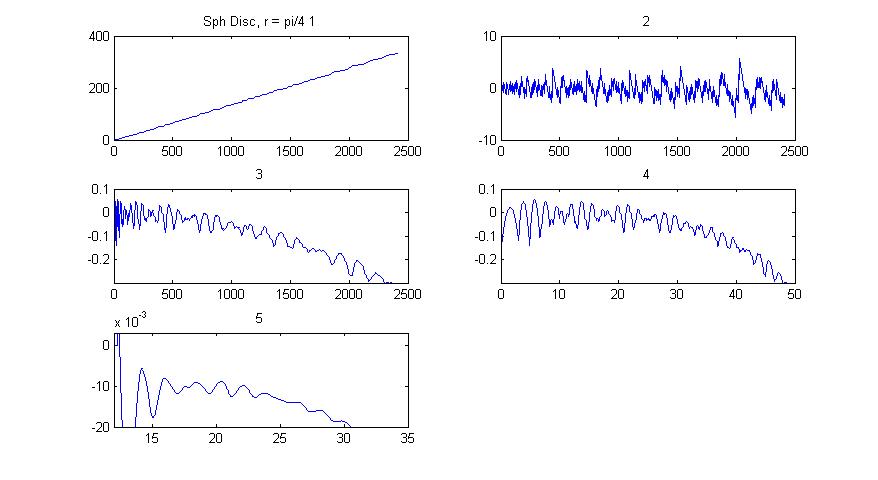}\\
\caption{Spherical disk with radius $r = \frac{\pi}{4}$}
\end{figure} 

%\begin{center}
%\includegraphics[scale=.40]{SphDisc_25DifferenceGraphs}\\
%Figure 5.8
%\end{center} 

\section{Discussion}\label{sec_discus}

We have presented strong experimental evidence for the conjectures, although certain statements in the conjectures that refer to almost periodicity are not candidates for verification by examining a bottom segment of the spectrum. Indeed, since the conjectures concern asymptotic behavior, it could have happened that the bottom segment would not have given a clue to the ultimate asymptotics. As it turned out we were lucky, and the conjectured asymptotic statements kicked in very early in the game. We are inclined to believe that this is not just luck, but that there is some paradigm at work here, to the effect that qualitative asymptotic statements about spectra can be refined to quantitative error estimates that would imply "early onset." We invite the reader to speculate about this possibility.

The conjectures, first put forth in \cite{strichartz2016}, were based on examining a collection of examples for which it is possible to compute the spectrum exactly. All of these examples exhibit a high degree of symmetry. It is always risky to jump to conclusions based on symmetric examples to the general case. Having now provided experimental evidence for surfaces that are not symmetric, convex, or even simply connected, we have a much firmer platform to support the conjectures. It should be kept in mind that the role of conjecture in mathematics is not always to light the way to future truth so much as to stimulate research on interesting problems. Even conjectures that eventually turned out to be incorrect have played an important role in the development of mathematics.

We note that in \cite{strichartz2016}, Conjecture \ref{conj_curv_neg} was stated only for flat surfaces. Indeed, there are no examples of hyperbolic surfaces for which it is possible to compute the spectrum exactly. Nevertheless, there is a great deal of interest in the spectra of hyperbolic surfaces, so we are pleased that our experiments support extending the conjecture to hyperbolic surfaces.

A very interesting question, which has not yet been explored in the literature, is the behavior of the differences of consecutive eigenvalues. Note that there is an immediate difference between the nature of differences on very symmetric surfaces and "generic" surfaces. Indeed, if the surface has a nonabelian symmetry group, then there will be eigenvalues with multiplicities greater than 1, so zero will be a difference that occurs often. In the generic case we do not expect any multiplicities greater than one. So we cannot expect a general answer that applies to all examples. We have gathered data on the differences for all our examples, and this may be found on the website \cite{murray2015results}, so it is possible that we have not been able to compute enough eigenvalues with enough accuracy to make a general pattern clear. We invite the reader to consider Figures 26-34 while thinking about this challenging question. In these figures, the $x$-axis is the distance between successive eigenvalues and the curve in the first figure if the number of differences less than or equal to $x$, while the second figure is a histogram of the successive differences. Note that eigenvalue spectra containing several values with multiplicities greater than 1 have the first curve start with a high proportion of the differences at 0 (in particular see Figure 33).

\begin{figure}[ht]
\centering
\includegraphics[width=\scalorDif\textwidth]{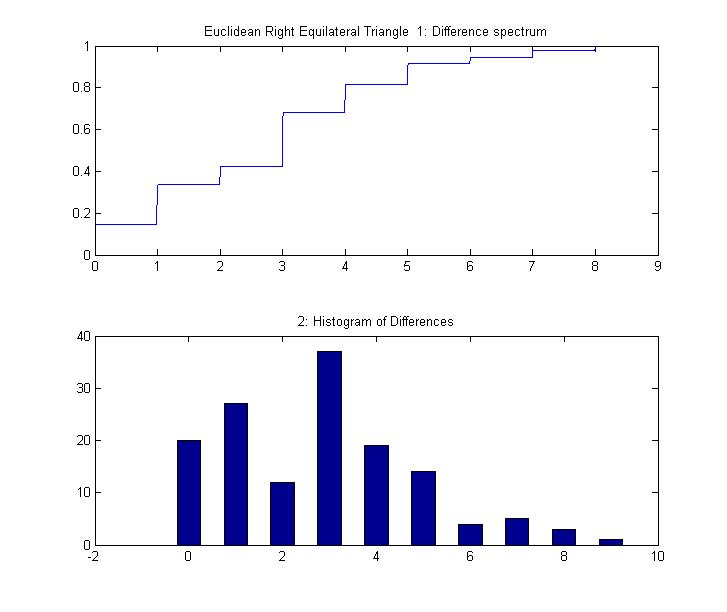}\\
\caption{Euclidean Right Isosceles Triangle with Dirichlet boundary conditions}
\end{figure}

\begin{figure}[ht]
\centering
\includegraphics[width=\scalorDif\textwidth]{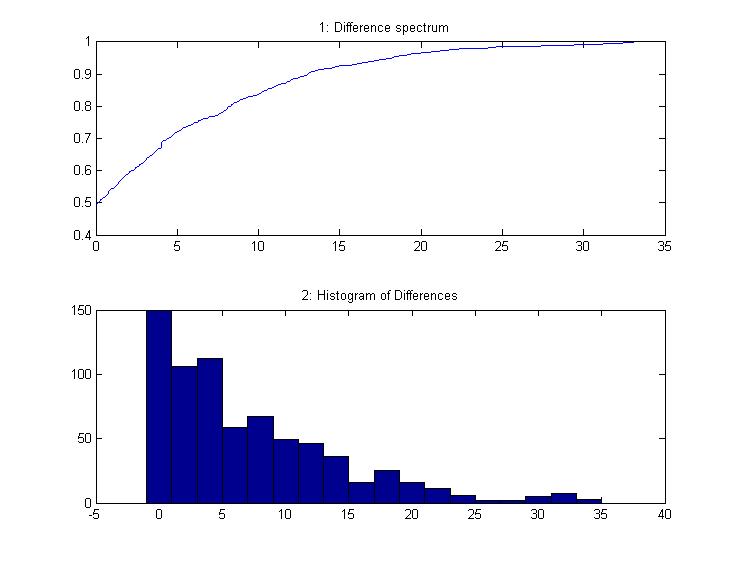}\\
\caption{Euclidean Disc with Dirichlet boundary conditions}
\end{figure} 

\begin{figure}[ht]
\centering
\includegraphics[width=\scalorDif\textwidth]{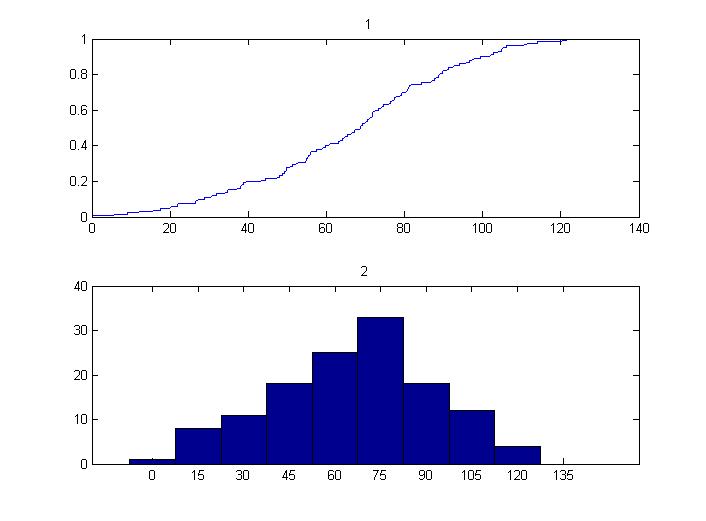}\\
\caption{General Euclidean Triangle with Dirichlet boundary conditions}
\end{figure} 

\begin{figure}[ht]
\centering
\includegraphics[width=\scalorDif\textwidth]{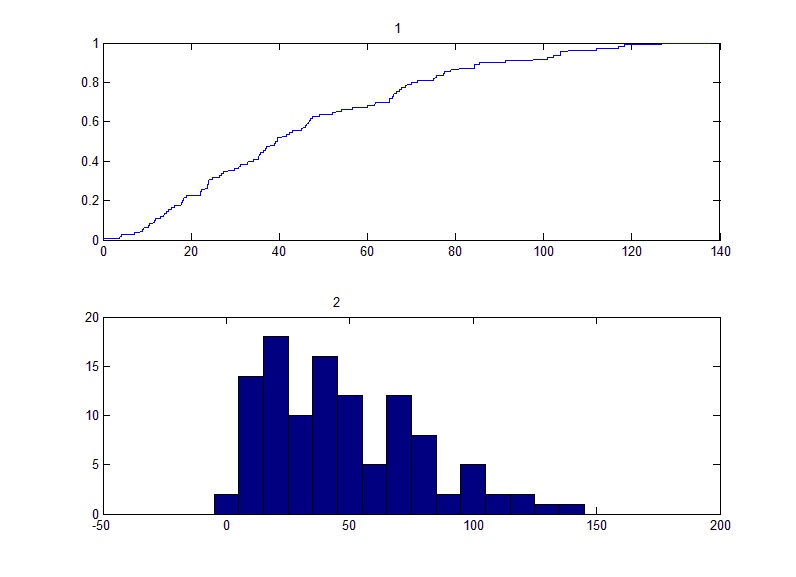}\\
\caption{Euclidean Region between Triangles with Dirichlet boundary conditions}
\end{figure} 

\begin{figure}[ht]
\centering
\includegraphics[width=\scalorDif\textwidth]{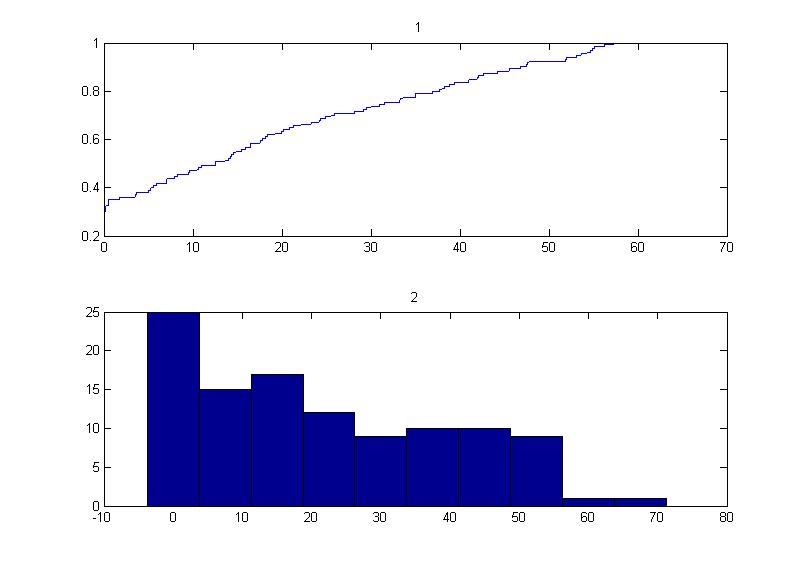}\\
\caption{Hyperbolic Equilateral Triangle ($\theta =\frac{\pi}{4}$) with Dirichlet boundary conditions}
\end{figure}

\begin{figure}[ht]
\centering
\includegraphics[width=\scalorDif\textwidth]{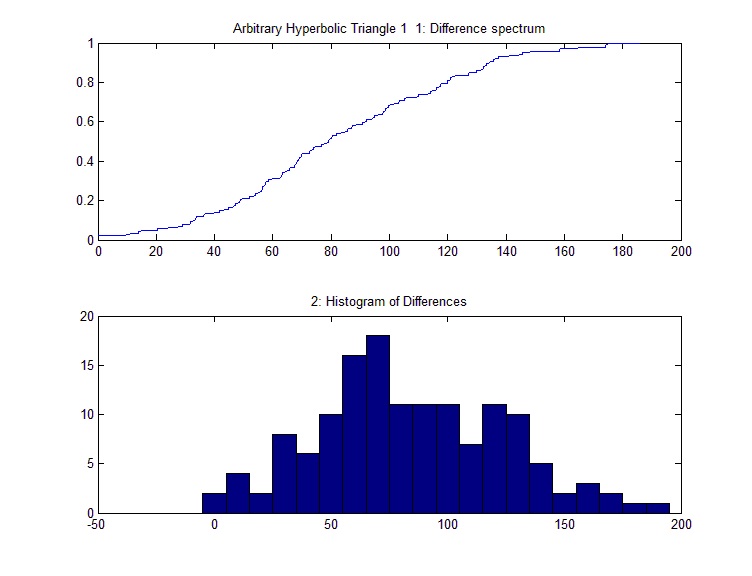}\\
\caption{General Hyperbolic Triangle with Dirichlet boundary conditions}
\end{figure} 

\begin{figure}[ht]
\centering
\includegraphics[width=\scalorDif\textwidth]{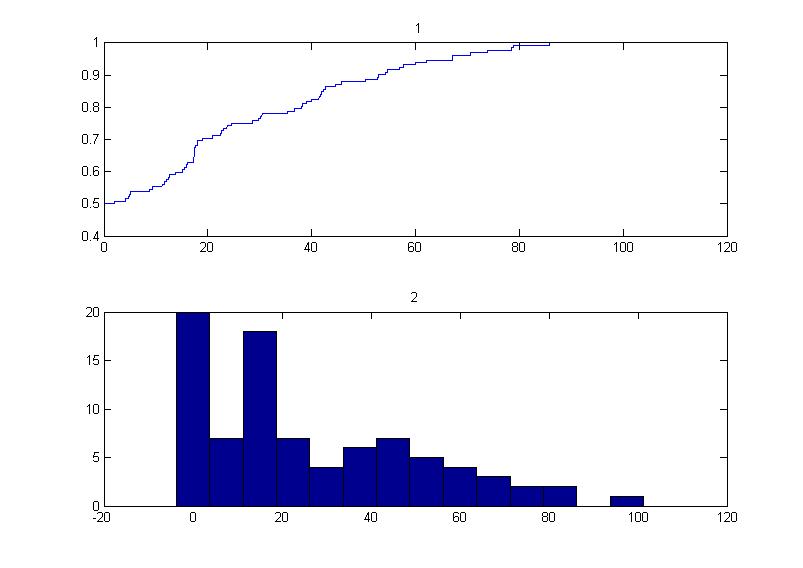}\\
\caption{Hyperbolic Disc ($R = \frac{1}{2}$) with Dirichlet boundary conditions}
\end{figure}

\begin{figure}[ht]
\centering
\includegraphics[width=\scalorDif\textwidth]{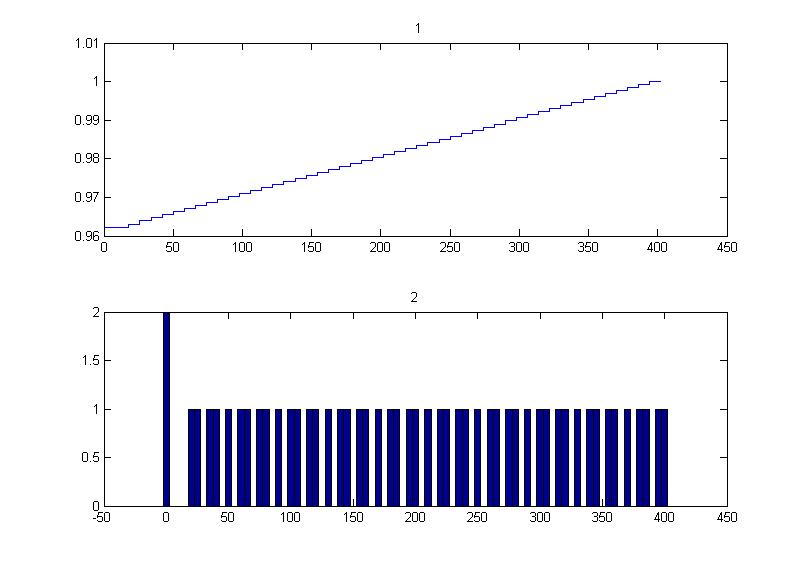}\\
\caption{Spherical Right Equilateral Triangle with Dirichlet boundary conditions}
\end{figure}  

\begin{figure}[ht]
\centering
\includegraphics[width=\scalorDif\textwidth]{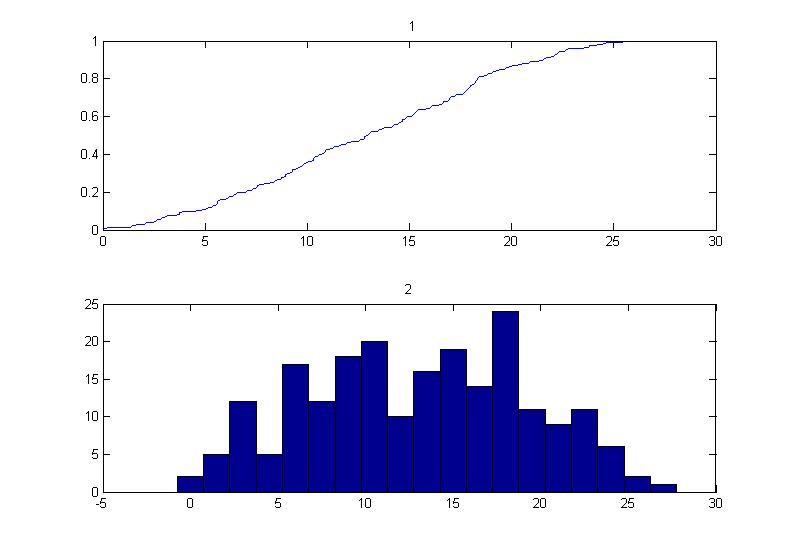}\\
\caption{Spherical General Triangle with Dirichlet boundary conditions}
\end{figure}  

\clearpage

\bibliographystyle{AIMS}
\bibliography{Main.bbl}

\medskip
% The data information below will be filled by AIMS editorial staff
Received xxxx 20xx; revised xxxx 20xx.
\medskip

\end{document}